\documentclass[12pt,a4paper]{article}

\usepackage{amssymb,amsmath,amsfonts,enumerate}
\usepackage{graphicx}

\usepackage[utf8]{inputenc}
\usepackage[english]{babel}

\numberwithin{equation}{section}
\newtheorem{theorem}{Theorem}[section]
\newtheorem{lemma}[theorem]{Lemma}%[sectionwp]
\newcommand{\MestnikovaIm}{\mathrm{\textsl{Im}\,}}
\newcommand{\MestnikovaRe}{\mathrm{\textsl{Re}\,}}

\title{Steady free surface potential flow of an ideal fluid\\
due to a singular sink on the flat bottom}
\author{A.A.~Mestnikova and V.N.~Starovoitov\\[3mm]
{\normalsize Lavrentyev Institute of Hydrodynamics, Novosibirsk, Russia}\\
{\small E-mail:\quad mestnikova@hydro.nsc.ru;\quad starovoitov@hydro.nsc.ru}}
\date{}

\begin{document}

\maketitle

\abstract{A two-dimensional steady problem of a potential free-surface flow of an ideal
incompressible fluid caused by a singular sink is considered. The sink is placed at the horizontal
bottom of the fluid layer. With the help of the Levi-Civita technique, the problem is rewritten
as an operator equation in a Hilbert space. It is proven that there exists a unique solution
of the problem provided that the Froude number is greater than some particular value.
The free boundary corresponding to this solution is investigated. It has a cusp over the sink and decreases
monotonically when going from infinity to the sink point. The free boundary is an analytic curve
everywhere except at the cusp point. It is established that the inclination angle of the free boundary
is less than $\pi/2$ everywhere except at the cusp point, where this angle is equal to $\pi/2$.
The asymptotics of the free boundary near the cusp point is investigated.}


\section{Introduction}\label{Mestnikova-sec1}
%\subsection{General statement of the problem}\label{Mestnikova-sec1.1}
In this paper, we investigate a two-dimensional steady problem of the
free surface potential flow of an incompressible perfect fluid over
the flat horizontal bottom. The flow is caused by a located at the bottom singular point sink
of the strength $m>0$. Let us introduce the rectangular
Cartesian coordinate system $(x,y)$ such that the bottom coincides with
the $x$-axis and the sink is located at its origin $O$.
We denote by $a_x$ and $a_y$ the corresponding components of a vector field
$\boldsymbol a=(a_x,a_y)$. It is assumed that the gravity force $\varrho\,\boldsymbol g$
acts on the fluid, where $\varrho=const$ is the density of the fluid,
$\boldsymbol g=(0,-g)$, and $g>0$ is the acceleration due to gravity. The problem is to find
the free surface $\varGamma$ and the velocity vector field $\boldsymbol v=(v_x,v_y)$.
If there is no sink ($m=0$), then the fluid is motionless and occupies
the horizontal layer of the depth $h$. In this case,
$\varGamma=\{(x,y)\in \mathbb{R}^2\;|\; y=h\}$.

\begin{figure}[h]
\begin{center}
\includegraphics[width=0.65\textwidth]{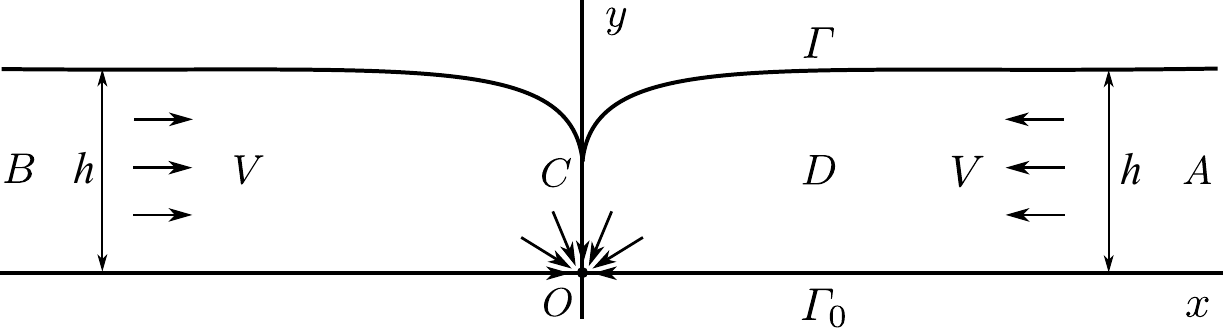}
\caption{\small General flow pattern in the physical plane.}
\label{Mestnikova-pic1}
\end{center}
\end{figure}

Let $D\subset \mathbb{R}^2$ be the domain occupied by the fluid whose
boundary consists of the free surface $\varGamma$ and
the bottom $\varGamma_0$ (see Fig.~\ref{Mestnikova-pic1}).
The velocity field satisfies in $D$ the steady Euler equations:
\begin{align}
\label{Mestnikova-1.1}
& v_x\partial_x\boldsymbol v+v_y\partial_y\boldsymbol v=-\frac{1}{\varrho}\,\nabla p +\boldsymbol g,
\\
\label{Mestnikova-1.2}
& \partial_x v_x+\partial_y v_y=0,
\end{align}
where $p$ is the pressure. We assume that the velocity field $\boldsymbol v$ is potential.
The corresponding equations are written in Section~\ref{Mestnikova-sec2.1}.
The Euler equations should be supplemented by boundary
conditions. Everywhere on the bottom $\varGamma_0$ except at the point $O$,
the no-flux condition is satisfied:
\begin{equation}\label{Mestnikova-1.3}
\varGamma_0\setminus \{O\}:\quad v_y=0.
\end{equation}
At the point $O$, the singular point sink of the strength $m$ is located, so
\begin{equation}\label{Mestnikova-1.4}
\Big(\boldsymbol v(x,y) + \frac{m}{2\pi}\,\frac{\boldsymbol r}{|\boldsymbol r|^2}\Big)\to \boldsymbol 0\quad\text{as}
\quad |\boldsymbol r|\to 0,
\end{equation}
where $\boldsymbol r=(x,y)$ and $|\boldsymbol r|^2=x^2+y^2$.

The upper boundary $\varGamma$ is unknown and, for this reason, two conditions
are imposed there. The first one is the kinematic condition:
\begin{equation}\label{Mestnikova-1.5}
\varGamma:\quad \boldsymbol v\cdot\boldsymbol n=0,
\end{equation}
where $\boldsymbol n$ is the normal vector to $\varGamma$. The second condition
assumes that the pressure is constant on $\varGamma$ and is equal to the atmospheric pressure. Since the Euler equations include only the gradient
of the pressure, we can impose the following condition:
\begin{equation}\label{Mestnikova-1.6}
\varGamma:\quad p=0.
\end{equation}

It remains to define the behavior of the fluid at the infinity points
$A$ and $B$ (see Fig.~\ref{Mestnikova-pic1}).
Depending on the value of $m$, various situations are possible.
The sink generates a perturbation of the free boundary that propagates as
waves to the right and to the left of the point above the sink. In this case, the velocity of the fluid does not tend to a certain limit at infinity. If the strength of the sink is sufficiently large, the velocity of the fluid on the surface exceeds the velocity of the waves and therefore we have a uniform flow at infinity. In this paper, we consider the case of sufficiently large $m$. This is the so called supercritical case.

Suppose that the fluid flow tends at infinity to a uniform flow whose depth is
$h$ and the value of the velocity is a constant $V> 0$:
\begin{equation}\label{Mestnikova-1.7}
\boldsymbol v(x,y)\to (\mp V, 0)\quad\text{as}\quad x\to\pm\infty .
\end{equation}
The constant $V$ is not arbitrary. Since the fluid is incompressible, this
quantity depends on the sink strength:
\begin{equation}\label{Mestnikova-1.8}
2Vh=\frac{m}{2}.
\end{equation}
On the right-hand side, $m$ is divided by two since only half of the sink drains the fluid from the flow domain $D$.

There are quite a large number of works devoted to the investigation of this and similar problems
qualitatively, in various approximations, and numerically. Since it is impossible to mention all of them,
we confine ourselves to papers directly related to our work. It seems that there are two key properties of the flow under consideration. The first one consists in the fact that, for sufficiently large values of the Froude number (to be defined below), there are no waves going to infinity. This property is a consequence of the monotonicity of the free surface when
going from infinity to the sink point.  The monotonicity was obtained numerically already in the first papers on the subject (see, e.g.,
\cite{Mestnikova-VBK,Mestnikova-H,Mestnikova-TVB,Mestnikova-FH}).  Notice that, in the case of the source, the waves do exist (see \cite{Mestnikova-LMCC}). In general, the waves also occur in the flows with small Froude numbers, i.e., in the so called subcritical case. This is not the subject of the present paper and we refer to the monograph \cite{Mestnikova-M}, where one can find an extensive bibliography on the question.
The second and not so obvious property of the flow is that the free boundary has a cusp over the sink for the sufficiently large Froude number. This fact was discovered numerically in the already cited works \cite{Mestnikova-VBK,Mestnikova-H,Mestnikova-TVB,Mestnikova-FH}. Here, we should notice that, for small Froude numbers, the stagnation point can occur on the free boundary over the sink (see \cite{Mestnikova-MVB,Mestnikova-HF}). At this point, the velocity of the fluid is equal to zero. The presence of the stagnation point is typical for the problem with the sink on the sloping bottom (see \cite{Mestnikova-DH} and also \cite{Mestnikova-VBK,Mestnikova-H}).

Since the purpose of this paper is the mathematical investigation of the problem, only qualitative results from the works cited above can be of use to us.
A quite effective mathematical technique was developed for the problem of the Stokes surface waves of extreme form (see \cite{Mestnikova-KN,Mestnikova-PT,Mestnikova-F}). This technique is based on the transition to the Nekrasov equation that exactly describes the free boundary. The study of this equation is based on the theory of positive solutions of nonlinear integral equations and on various results of harmonic analysis.
We cannot apply this technique without changes because the Nekrasov equations for our case and for the case of periodic waves differ.
Another distinction is in the proof of the solvability.
In the theory of the Stokes waves, a continuous one-parameter family of solutions of the Nekrasov equation is constructed.
At a some value of the parameter, the solution is easy to find and the Stokes waves correspond to another value of the parameter.
In our case, there is no value of the parameter at which the proof of the solvability of the problem would be trivial. Of course, if the Froude number is zero, there is no sink and the zero velocity field together with the flat free boundary form the solution. However, this solution does not belong to a continuous branch of solutions.
Finally, the theory of the Stokes waves deals with bounded functions while we have the singularity of the velocity field due to the sink.

The paper is organized as follows. In the next section, we rewrite the problem in a more appropriate
for investigation form. As we deal with the two-dimensional problem, it is quite natural to use the complex variable.
With the help of conformal mappings and the generalized Levi-Civita method, we write down an equation of the Nekrasov type on the unit circle.
The complete equivalent formulation of the problem is given in Section~\ref{Mestnikova-sec2.3}. It consists of two equations one of which
is the Nekrasov type integro-differential equation and the second is the integral equation that is one of the Hilbert inversion formulas.
In Section~\ref{Mestnikova-sec3}, we formulate the problem as one nonlinear operator equation in the function space $L^2(0,\pi/2)$ and prove
its unique solvability (Theorems~\ref{Mestnikova-t3.7} and \ref{Mestnikova-t3.8}). In Section~\ref{Mestnikova-sec4}, the differentiability (Theorem~\ref{Mestnikova-t4.3}) and
even the analyticity (Theorem~\ref{Mestnikova-t4.5}) of the free boundary is established. Besides, it is proven there that the inclination angle
of the free boundary is less than $\pi/2$ everywhere except at the point over the sink, where this angle is equal to $\pi/2$, i.e., the boundary
has the cusp (Theorem~\ref{Mestnikova-t4.7}). The asymptotics of the free boundary near the cusp point is investigated (Theorem~\ref{Mestnikova-t4.8} and Theorem~\ref{Mestnikova-t4.9}).
All these results are obtained for the so-called supercritical case, when the Froude number exceeds a certain positive value.
Finally, in Appendix, we give a proof of a representation for the kernel of an integral operator similar to the Hilbert transform.
This representation is used in many papers, but we have found its proof only in \cite{Mestnikova-N1} which is a fairly rare book.
Notice that our proof differs from that in \cite{Mestnikova-N1}.

\section{Statement of the problem on the unit circle}\label{Mestnikova-sec2}

The statement of the problem given in the previous section is quite difficult to deal with.
It is unclear how to carry out even a numerical implementation. For this reason, in this section, we rewrite the problem
in a more appropriate form.

\subsection{Complex variable formulation of the problem}\label{Mestnikova-sec2.1}

\noindent
\emph{Equations of the potential flow in $D$.}\\
As already mentioned above, we suppose that the flow is potential.
Thus, there exists a scalar function $\varphi=\varphi(x,y)$
such that $\boldsymbol v =\nabla\varphi$, i.e.,
$v_x=\partial_x\varphi$ and  $v_y=\partial_y\varphi$.
Due to equation \eqref{Mestnikova-1.2}, the potential $\varphi$ is a harmonic function:
\begin{equation}\label{Mestnikova-2.1}
\Delta\varphi =0\quad\text{in}\quad D.
\end{equation}
Besides that, equation \eqref{Mestnikova-1.2} implies that there exists a so
called stream function $\psi=\psi(x,y)$, such that $v_x=\partial_y\psi$ and $v_y=-\partial_x\psi$. This function is also harmonic:
\begin{equation}\label{Mestnikova-2.2}
\Delta\psi =0\quad\text{in}\quad D.
\end{equation}
The potential and the stream function are defined up to an arbitrary
additive constant.

If we have found the potential or the stream function, then we are
able to calculate the velocity vector field and, after that, the pressure can be found from equation \eqref{Mestnikova-1.1}. The pressure, however, does not interest us, and we will not define it. Therefore, it suffices to restrict ourselves to  solving equation \eqref{Mestnikova-2.1} or \eqref{Mestnikova-2.2} with appropriately rewritten boundary conditions. Here, we encounter a difficulty due to the fact that the boundary condition \eqref{Mestnikova-1.6} includes the pressure and must be written in a different form.

\medskip\noindent
\emph{Boundary conditions at $\varGamma_0$ and $\varGamma$.}\\
Streamline is a curve along which the stream function is constant.
It is not difficult to see that streamline can be also defined as a curve
tangent to the velocity vector. As follows from the boundary conditions \eqref{Mestnikova-1.3} and \eqref{Mestnikova-1.5}, the boundaries $\varGamma_0$ and $\varGamma$
are streamlines. Notice, however, that each of these boundaries consists of two streamlines. Let us introduce the following notations:
\begin{gather*}
\varGamma_0^+=\{(x,y)\in\varGamma_0\;|\; x>0\},
\quad\varGamma_0^-=\{(x,y)\in\varGamma_0\;|\; x<0\},
\\
\varGamma^+=\{(x,y)\in\varGamma\;|\; x>0\},
\quad\varGamma^-=\{(x,y)\in\varGamma\;|\; x<0\}.
\end{gather*}
Since the problem is symmetric with respect to the $y$-axis, the velocity
vector is directed vertically (downwards) on the interval $(O,C)$ (see Fig.~\ref{Mestnikova-pic1}), i.e., it is tangent to this axis. Therefore, this interval is also a streamline.
The function $\psi$ is defined up to an additive constant,
therefore, without loss of generality, we can assume that it vanishes at the point $C$ and, as a consequence, that
\begin{equation}\label{Mestnikova-2.3}
\psi|_{\varGamma}=\psi|_{(O,C)}=0.
\end{equation}

Let us consider the streamlines $\varGamma_0^+$ and $\varGamma_0^-$.
The difference of the stream function values at two points represents
the flow rate of the fluid through any curve connecting these points. Therefore,
\begin{equation}\label{Mestnikova-2.4}
\psi|_{\varGamma_0^+}=Vh,\quad \psi|_{\varGamma_0^-}=-Vh.
\end{equation}

It remains to rewrite the boundary condition \eqref{Mestnikova-1.6} in a form that does not include the pressure. To this end, we use the Bernoulli principle which states that the quantity $|\boldsymbol v|^2/2 + p/\varrho + gy$ is constant along any streamline. Since $p = 0$ on the streamlines $\varGamma^\pm$, we get the following relation:
\begin{equation}\label{Mestnikova-2.5}
\frac{1}{2}\,|\boldsymbol v(x,y)|^2 + gy=K=const\quad\text{as}\quad (x,y)\in\varGamma.
\end{equation}
Due to \eqref{Mestnikova-1.7}, $K=V^2/2 + gh$.

\medskip\noindent
\emph{Dimensionless formulation of the problem.}\\
Despite the fact that the problem contains several dimensional parameters
such as the strength of the sink $m$, the depth of the undisturbed fluid layer
$h$, the gravitational acceleration $g$, the fluid flow is determined by one dimensionless parameter, the Froude number
$$
\mathrm{Fr}=\frac{V}{\sqrt{gh}}.
$$
The quantity $V$ is not an original parameter of the problem. By using \eqref{Mestnikova-1.8}, $\mathrm{Fr}$ can be expressed in terms of
$m$, $h$, and $g$:
$$
\mathrm{Fr}^2=\frac{m^2}{16gh^3}.
$$
For brevity, we use in the paper the reduced Froude number $\alpha=\mathrm{Fr}^2/2$.

Let us take $ h $ and $ V $ as the characteristic units of the length and the velocity, respectively, and leave the previous notations for all dimensionless quantities. The dimensionless potential $\varphi$ and stream function $\psi$ satisfy \eqref{Mestnikova-2.1}, \eqref{Mestnikova-2.2}, and \eqref{Mestnikova-2.3} with dimensionless $x$ and $y$. Similarly, for the dimensionless velocity, we have the following expressions: $v_x=\partial_x\varphi=\partial_y\psi$ and $v_y=\partial_y\varphi=-\partial_x\psi$.
Conditions \eqref{Mestnikova-2.4}, \eqref{Mestnikova-2.5}, \eqref{Mestnikova-1.4}, and \eqref{Mestnikova-1.7}
take the form:
\begin{equation*}
\psi|_{\varGamma_0^+}=1,\quad \psi|_{\varGamma_0^-}=-1,
\end{equation*}
\begin{equation}\label{Mestnikova-2.6}
\alpha\,|\boldsymbol v(x,y)|^2 + y=\alpha+1\quad\text{for}\quad (x,y)\in\varGamma,
\end{equation}
\begin{equation}\label{Mestnikova-2.7}
\Big(\boldsymbol v(x,y) + \frac{2}{\pi}\,\frac{\boldsymbol r}{|\boldsymbol r|^2}\Big)\to \boldsymbol 0\quad\text{as}
\quad |\boldsymbol r|\to 0,
\end{equation}
\begin{equation}\label{Mestnikova-2.8}
\boldsymbol v(x,y)\to (\mp 1, 0)\quad\text{as}\quad x\to\pm\infty .
\end{equation}
Notice also that $y\to 1$ as $x\to\pm\infty$ for $(x,y)\in\varGamma$.

\medskip\noindent
\emph{Complex variable formulation of the problem.}\\
Since we are solving a two-dimensional stationary problem with harmonic functions,
it will be convenient to employ the complex variable functions theory. The functions $\varphi$ and $\psi$ satisfy in $D$ the Cauchy --- Riemann equations:
$$
\partial_x \varphi=\partial_y \psi, \quad \partial_y \varphi=-\partial_x \psi,
$$
therefore, the complex function $w=\varphi + i\psi$  is a
holomorphic function of the complex variable $z=x+iy$ in the domain $D$.
This function is called the complex potential.
The complex function $v=v_x-iv_y$ is the derivative of the complex potential
with respect to $z$: $v=d w/d z$. This function is called the complex velocity.
The complex variable formulation of the Bernoulli principle \eqref{Mestnikova-2.6} looks as follows:
\begin{equation}\label{Mestnikova-2.9}
\alpha\,|v(z)|^2 + \MestnikovaIm z=\alpha+1 \quad\text{for}\quad z\in\varGamma.
\end{equation}
Condition \eqref{Mestnikova-2.7} takes the form:
\begin{equation*}
\Big(v(z) + \frac{2}{\pi}\,\frac{1}{z}\Big)\to 0\quad\text{as}
\quad z\to 0.
\end{equation*}
The complex formulation of the remaining boundary conditions at $\varGamma_0$ and $\varGamma$ does not cause difficulties.

\subsection{The formulation of the problem on the unit circle}\label{Mestnikova-sec2.2}

\noindent
\emph{The problem in a half-disk.}\\
Let us denote by $F$ the conformal mapping of the flow domain $D$ to
the upper half $D^*$ of the unit disk $D_\circ^*$ centered at the origin $O^*$ of the plane
of the complex variable $t$. We require additionally that $F$ maps the point $O$ to the point $O^*(t = 0)$, the point $C$ to the point $C^*(t = i)$,
the segment $[O,C]$ to the segment $[O^*,C^*]$,
the points $A$ and $B$ at infinity to the points $A^*(t = 1)$ and $B^*(t = -1)$, respectively. Due to the symmetry of the problem with respect to the
$y$-axis, such a conformal mapping $F$ exists and $\MestnikovaRe F(z_1)=-\MestnikovaRe F(z_2)$
for all points $z_1$ and $z_2$ such that $\MestnikovaRe z_1=-\MestnikovaRe z_2$.
The images of the free boundary $\varGamma$ and the bottom $\varGamma_0$ under the mapping $F$ are the upper half $\varGamma^*$ of the unit circle and the horizontal diameter $\varGamma_0^*=[B^*,A^*]$, respectively.
Since the domain $D$ is unknown, we can not determine the mapping $F$ at this stage of solving the problem. However, for each domain $D$, such a mapping exists and is uniquely determined.

\begin{figure}[h]
\begin{center}
\includegraphics[width=0.45\textwidth]{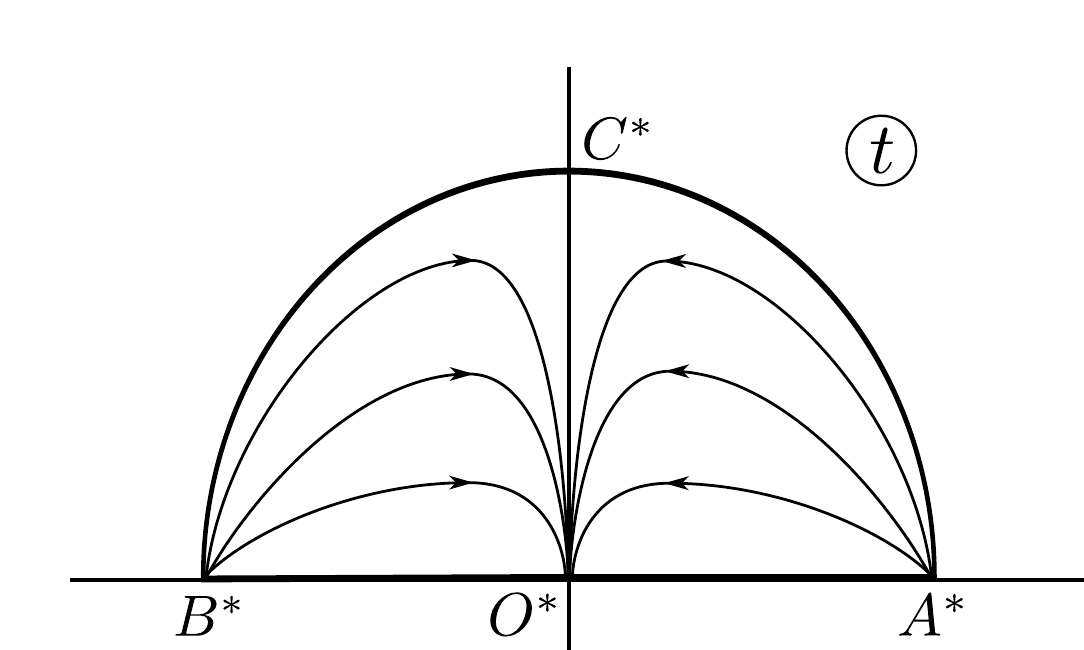}
\caption{\small The flow pattern in the half-disk $D^*$.}
\label{Mestnikova-pic2}
\end{center}
\end{figure}

Denote by $G$ the inverse mapping to $F$. The function $w^*(t) = w\big(G(t)\big)$ is the complex potential of the flow in the domain $D^*$.
The no-flux conditions are satisfied on the boundaries
$\varGamma^*$ and $\varGamma_0^*$. There is a point sink at the point $O^*$ and  identical point sources at the points $A^*$ and $B^*$ (see Fig.~\ref{Mestnikova-pic2}).
Notice that $w(z)=w^*\big(F(z)\big)$ is the complex potential of the flow
in $D$ provided that $w^*$ is the complex potential in $D^*$.
It is not difficult to find that
\begin{equation}\label{Mestnikova-2.10}
w^*(t)=\frac{2}{\pi}\,\log\Big(\frac{t^2-1}{2t}\Big)-i
=\frac{2}{\pi}\,\big(\log(t+1)+\log(t-1)-\log 2t-i\pi/2\big) .
\end{equation}

Thus, the flow in the half-disk $D^*$ is completely defined.
If we knew the mapping $F$, then we could determine the complex potential of the flow in the plane of the variable $z$ and, thus, solve the problem.  On the other hand, this mapping is uniquely determined by the domain $D$ or equivalently by
its upper boundary $\varGamma$ which is a priori unknown.
To determine $\varGamma$, we use the Bernoulli equation \eqref{Mestnikova-2.6} and the Levi-Civita approach which
is described below.

\medskip\noindent
\emph{The Levi-Civita approach.}\\
Since $w(z)=w^*\big(F(z)\big)$,
$$
v(z)=\frac{dw}{dz}(z)=\frac{dw^*}{dt}\big(F(z)\big)\,\frac{dF}{dz}(z).
$$
Therefore,
\begin{equation*}
v\big(G(t)\big)=\frac{dw^*}{dt}(t)\,\frac{dF}{dz}\big(G(t)\big)=
\frac{dw^*}{dt}(t)\,\Big(\frac{dG(t)}{dt}\Big)^{-1}.
\end{equation*}
This equality can be also rewritten as follows:
\begin{equation}\label{Mestnikova-2.11}
\frac{dG(t)}{dt}=\frac{1}{v\big(G(t)\big)}\,\frac{dw^*}{dt}(t).
\end{equation}
If we knew the function $v\big(G(t)\big)$, then we could find $G$ from this equation.  The function $v\big(G(t)\big)$ is holomorphic in $D^*$ and has a pole of the first order at the point $t=0$. Besides that, it tends to $-1$ and $1$
as $t$ tends to $A^*(t=1)$ and $B^*(t=-1)$, respectively.
Due to the symmetry of the problem, this function is purely imaginary
at the imaginary axis: $\MestnikovaRe v\big(G(t)\big)=0$ for $t\in(0,i]$.

Everywhere on the diameter $A^*O^*B^*$, except at the point $O^*(t = 0)$, the no-flux condition is fulfilled: $\MestnikovaIm v\big(G(t)\big) = 0$. Therefore, using the Schwartz symmetry principle, we can analytically continue the function $v\big(G(t)\big)$ to the lower half of the unit disk $D_\circ^*$. We denote the continued function by $u(t)$. The function $u(t)$ has the pole of the first order at the point $t=0$.

Following the Levi-Civita approach, we introduce the function
$$
\varOmega(t)=\hat{\theta}(t)+i\hat{\tau}(t)
$$
with real $\hat{\theta}$ and $\hat{\tau}$ such that
\begin{equation}\label{Mestnikova-2.12}
u(t)=-\frac{2}{\pi}\,\frac{1}{t}\, e^{-i\varOmega(t)}.
\end{equation}
It is not difficult to see that $\varOmega(t)$ is holomorphic in $D_\circ^*=\{t\;|\; |t|<1\}$.
Notice that the function $\hat{\theta}$ is defined up to an additive
constant $2\pi k$, where $k$ is an integer number.
Further, in \eqref{Mestnikova-2.17}, we fix one branch of this function by choosing its value at the point $t=1$.

\medskip\noindent
\emph{The Bernoulli equation on $\varGamma^*$.}\\
The Bernoulli identity \eqref{Mestnikova-2.9} has the following form on $\varGamma^*$:
\begin{equation}\label{Mestnikova-2.13}
\alpha\,|v(G(t))|^{2}+ \MestnikovaIm G(t)=\alpha+1.
\end{equation}
Since $\varGamma^*$ is a part of the unit circle, $t=e^{i\sigma}$ with $\sigma\in[0,\pi]$ for every $t\in \varGamma^*$.  Substitution of this representation into \eqref{Mestnikova-2.13} and differentiation with respect to $\sigma$ give:
\begin{equation}\label{Mestnikova-2.14}
\alpha\,\frac{d\big|v\big(G(e^{i\sigma})\big)|^{2}}{d\sigma}+
\MestnikovaIm\frac{dG(e^{i\sigma})}{d\sigma}=0
\quad\text{for}\quad \sigma\in (0,\pi).
\end{equation}
Let us calculate $\MestnikovaIm dG(e^{i\sigma})/d\sigma$. At first, we note that
$$
\frac{dG(e^{i\sigma})}{d\sigma}=\frac{dG}{dt}\Big|_{t=e^{i\sigma}}\,ie^{i\sigma}.
$$
Taking into account that $v(G(t))=u(t)$ and using equalities \eqref{Mestnikova-2.11}, \eqref{Mestnikova-2.12}, and \eqref{Mestnikova-2.10}, we find that
$$
\frac{dG(t)}{dt}=\frac{1}{u(t)}\,\frac{dw^*}{dt}=-\frac{\pi t}{2}
\, e^{i\varOmega(t)}\, \frac{2}{\pi t}\,\frac{t^{2}+1}{t^{2}-1}=
-e^{i\varOmega(t)}\,\frac{t^{2}+1}{t^{2}-1}.
$$
Therefore,
$$
\frac{dG(e^{i\sigma})}{d\sigma}=-  e^{i(\sigma+\varOmega(t))}\Big|_{t=e^{i\sigma}}\,\cot\sigma.
$$
We are interested in the values of the functions $\hat{\tau}$ and $\hat{\theta}$ at the points $t=e^{i\sigma}$. For this reason, we introduce the following
notations:
$$
\tau(\sigma)=\hat{\tau}(e^{i\sigma}),\quad \theta(\sigma)=\hat{\theta}(e^{i\sigma}).
$$
Thus,
$$
\frac{dG(e^{i\sigma})}{d\sigma}=
- e^{-\tau(\sigma)} \, e^{i(\sigma+\theta(\sigma))}\, \cot\sigma
$$
and, as a consequence,
$$
\MestnikovaIm\frac{dG(e^{i\sigma})}{d\sigma}=-  e^{-\tau(\sigma)} \, \sin(\sigma+\theta(\sigma))\,\cot\sigma.
$$

Let us consider now the first term on the right-hand side of \eqref{Mestnikova-2.14}.
As it follows from \eqref{Mestnikova-2.12},
$$
\big|v\big(G(e^{i\sigma})\big)\big|^2=\big|u(e^{i\sigma})\big|^2=
\frac{4}{\pi^2}\, e^{2\tau(\sigma)}.
$$
Substitution of the obtained expressions in \eqref{Mestnikova-2.14} leads to the equation
$$
\frac{4\alpha}{\pi^2}\,\frac{de^{2\tau(\sigma)}}{d\sigma}
- e^{-\tau(\sigma)} \, \sin(\sigma+\theta(\sigma))\, \cot\sigma=0
$$
and, as a consequence, to the following form of the Bernoulli equation:
\begin{equation}\label{Mestnikova-2.15}
e^{3\tau(\sigma)}\tau\,'(\sigma)-\frac{\pi^2}{8\alpha}\, \sin(\sigma+\theta(\sigma))\,\cot\sigma=0,
\end{equation}
where $\tau\,'(\sigma)=d\tau(\sigma)/d\sigma$.

\medskip\noindent
\emph{Values of the functions $\theta$ and $\tau$ at some points.}\\
Further, we need to know the values of the functions $\theta$ and $\tau$ at some points. We begin with  the investigation of the function $\hat{\theta}(t)$ on the diameter $[B^*,A^*]$, where $\MestnikovaIm t=0$ and $\MestnikovaIm u(t)=0$ everywhere except at the point $t=0$. Let $s$ be an arbitrary real number from the interval $(0,1]$. Then
$$
\MestnikovaIm u(t)\big|_{t=s}=\frac{2}{\pi}\,\frac{1}{s}\, e^{\hat{\tau}(s)}\,\sin\hat{\theta}(s)=0.
$$
This means that $\sin\theta(s)=0$ for $s\in(0,1]$. On the other hand, due to \eqref{Mestnikova-2.8}, $u(s)\to -1$ as $s\to 1$, that is
\begin{equation}\label{Mestnikova-2.16}
-\frac{2}{\pi}\,e^{\hat{\tau}(1)}\,\cos\hat{\theta}(1)=-1.
\end{equation}
Therefore, $\sin\hat{\theta}(1)=0$ and $\cos\hat{\theta}(1)>0$. We choose the following solution of this system:
\begin{equation}\label{Mestnikova-2.17}
\hat{\theta}(t)\big|_{t=1}=0.
\end{equation}
By this choice, we have fixed a branch of the function $\hat{\theta}$, as mentioned after equation \eqref{Mestnikova-2.12}.

Since the function $\hat{\theta}$ is continuous and $\sin\hat{\theta}(s)=0$ for $s\in(0,1]$, we deduce that $\hat{\theta}(t)=0$ for $t\in [O^*,A^*]$.
By the same arguments, we get that $\hat{\theta}(t)=0$ for $t\in [B^*,O^*]$.
Thus,
\begin{equation}\label{Mestnikova-2.18}
\hat{\theta}(t)=0\quad\text{as}\quad \MestnikovaIm t=0
\end{equation}
and, in particular,
\begin{equation}\label{Mestnikova-2.19}
\theta(0)=\theta(\pi)=0.
\end{equation}
Besides, due to \eqref{Mestnikova-2.16} and a similar relation for $t=-1$, we get
that
\begin{equation}\label{Mestnikova-2.20}
\tau(0)=\tau(\pi)=\log\frac{\pi}{2}.
\end{equation}

Let us now consider  $t=is$ with $s\in(-1,0)\cup(0,1)$.
As we have already found out, $\MestnikovaRe u(t) = 0 $
and, as an easy consequence, $\sin\hat{\theta}(t) = 0 $ for such values of $t$. Since $\hat{\theta}$ is a continuous function and $\hat{\theta}(0) = 0$,
\begin{equation*}
\hat{\theta}(t)=0\quad\text{as}\quad \MestnikovaRe t=0.
\end{equation*}
In particular,
\begin{equation*}
\theta(\pi/2)=\theta(3\pi/2)=0.
\end{equation*}

\medskip\noindent
\emph{The Nekrasov equation.}\\
At first, we note that equation \eqref{Mestnikova-2.15} can be written in the following
form:
\begin{equation*}
\frac{d}{d\sigma}e^{3\tau(\sigma)}=\frac{3\pi^2}{8\alpha}\, \sin(\sigma+\theta(\sigma))\,\cot\sigma.
\end{equation*}
Integration of this equation from $0$ to an arbitrary $\sigma$ together with
the condition \eqref{Mestnikova-2.20} gives
\begin{equation*}
e^{3\tau(\sigma)}=\frac{\pi^3}{8}+\frac{3\pi^2}{8\alpha}\,\int_0^\sigma
\sin(s+\theta(s))\,\cot s \, ds.
\end{equation*}
Substitution of this expression into \eqref{Mestnikova-2.15} leads to the following
Nekrasov type equation:
\begin{equation}\label{Mestnikova-2.21}
\tau\,'(\sigma)=
\frac{\sin(\sigma+\theta(\sigma))\,\cot\sigma}{\alpha\pi+3\int_{0}^\sigma
\sin(s+\theta(s))\,\cot s\, ds}.
\end{equation}
A similar equation for the surface waves was derived by A.I.~Nekrasov (see \cite{Mestnikova-N2} and \cite{Mestnikova-N1}).

\medskip\noindent
\emph{The Hilbert inversion formulas on the quarter of the circle.}\\
Equation \eqref{Mestnikova-2.21} includes two unknown functions, so it is necessary to find some other relation for them. We use the fact that these functions are traces on the unit circle of the real and imaginary parts of the holomorphic in the unit
disk $D_\circ^*$ function $\varOmega$. This means that the Hilbert inversion formulas hold for the functions $\tau$ and $\theta$:
\begin{align}
\notag
& \tau(\sigma)=-\frac{1}{2\pi}\,\int_0^{2\pi} \theta(s)\,
\cot\frac{s-\sigma}{2}\, ds + \hat{\tau}_0,
\\
\label{Mestnikova-2.22}
& \theta(\sigma)=\frac{1}{2\pi}\,\int_0^{2\pi} \tau(s)\,
\cot\frac{s-\sigma}{2}\, ds + \hat{\theta}_0,
\end{align}
where $\sigma\in[0,2\pi)$, $\hat{\tau}_0=\hat{\tau}|_{z=0}$ and $\hat{\theta}_0=\hat{\theta}|_{z=0}$.
Here, the integrals are understood in the sense of
the Cauchy principal value.

Notice that the constant $\hat{\tau}_0$ is unknown. At the same time, as it follows from \eqref{Mestnikova-2.18}, $\hat{\theta}_0=0$ and equation \eqref{Mestnikova-2.22} is completely defined. Since we need only one equation in addition to \eqref{Mestnikova-2.21}, we use \eqref{Mestnikova-2.22} which in our case looks as follows:
\begin{equation}\label{Mestnikova-2.23}
\theta(\sigma)=\frac{1}{2\pi}\,\int_0^{2\pi} \tau(s)\,
\cot\frac{s-\sigma}{2}\, ds.
\end{equation}

Our next goal is to rewrite this formula in such a way that it
will include only the integral over the interval $[0,\pi/2]$.
To this end, we employ the symmetry properties of the functions $\tau$ and $\theta$ that follow from the symmetry of the problem:
\begin{subequations}
\label{Mestnikova-2.24}
\begin{gather}
\label{Mestnikova-2.24a}
\tau(\pi/2+\sigma)=\tau(\pi/2-\sigma),\quad
\theta(\pi/2+\sigma)=-\theta(\pi/2-\sigma),
\\
\tau(\sigma)=\tau(-\sigma),\quad
\theta(\sigma)=-\theta(-\sigma),
\\
\label{Mestnikova-2.24c}
\tau(\pi+\sigma)=\tau(\sigma),\quad
\theta(\pi+\sigma)=\theta(\sigma)
\end{gather}
\end{subequations}
for all $\sigma\in[0,2\pi]$.

Due to \eqref{Mestnikova-2.24c}, by making the change of variable $s=\pi+\gamma$,
we find that
$$
\int\limits_\pi^{2\pi} \tau(s)\cot\frac{s-\sigma}{2} \,ds=
\int\limits_0^\pi \tau(\pi+\gamma)\,\cot\Big(\frac{\pi}{2}+\frac{\gamma-\sigma}{2}\Big)\,d\gamma=
-\int\limits_0^\pi \tau(\gamma)\,\tan\frac{\gamma-\sigma}{2}\, d\gamma.
$$
Since $\cot\alpha-\tan\alpha=2\cot 2\alpha$, equation \eqref{Mestnikova-2.23} takes the form:
\begin{equation*}
\theta(\sigma)=\frac{1}{\pi}\int_0^\pi \tau(s)\cot(s-\sigma)\,ds.
\end{equation*}

Since $\tau(\pi-\sigma)=\tau(\sigma)$, we easily find that
$$
\int_{\pi/2}^\pi \tau(s)\cot(s-\sigma)\,ds=-\int_{0}^{\pi/2} \tau(s)\cot(s+\sigma)\,ds.
$$
Therefore,
\begin{equation*}
\theta(\sigma)=\frac{1}{\pi}\int_0^{\pi/2} \tau(s)\,\big(\cot(s-\sigma)-\cot(s+\sigma)\big)\,ds.
\end{equation*}

\subsection{Complete formulation of the original problem in terms of the functions $\tau$ and $\theta$}\label{Mestnikova-sec2.3}
As an intermediate result, we will gather all the equations obtained and show how one can find the solution of the original problem. Recall that the original problem was to determine the free boundary $\varGamma$ of the flow domain $D$ and the velocity field of this flow. We construct the solution in several steps.

\medskip\noindent
\emph{Step~1.}
The first step is to find functions $\tau$ and $\theta$ that satisfy
the following equations:
\begin{equation}\label{Mestnikova-2.25}
\tau\,'(\sigma)=
\frac{\sin(\sigma+\theta(\sigma))\,\cot\sigma}{\alpha\pi+3\int_{0}^\sigma
\sin(s+\theta(s))\,\cot s\, ds},
\end{equation}
\begin{equation}\label{Mestnikova-2.26}
\theta(\sigma)=\frac{1}{\pi}\int_0^{\pi/2} \tau(s)\,\big(\cot(s-\sigma)-\cot(s+\sigma)\big)\,ds
\end{equation}
for $\sigma\in (0,\pi/2)$. The first of them is the Nekrasov type equation and the second is a form of the Hilbert inversion formula.
Besides that, these functions satisfy the following boundary conditions:
\begin{equation}\label{Mestnikova-2.27}
\tau(0)=\log\frac{\pi}{2},\quad
\theta(0)=\theta(\pi/2)=0.
\end{equation}
Notice that the boundary conditions for $\theta$ are taken from \eqref{Mestnikova-2.19}, however, they follow also
from \eqref{Mestnikova-2.26}.
After we have found the functions $\tau$ and $\theta$ on the interval
$[0,\pi/2]$, they can be determined on $[0,2\pi]$ with the help of the symmetry
properties \eqref{Mestnikova-2.24}.

\begin{figure}[h]
\begin{center}
\includegraphics[width=0.55\textwidth]{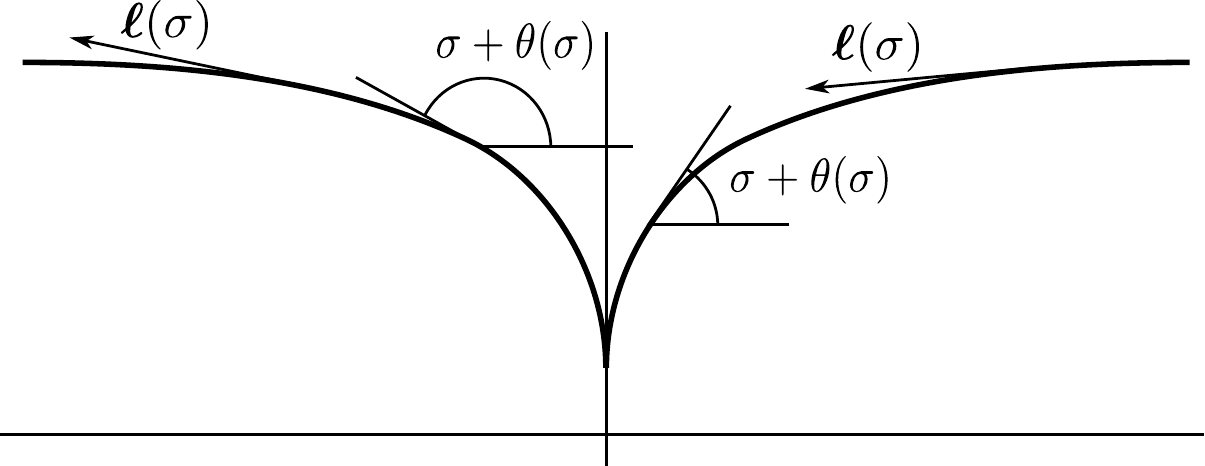}
\caption{\small The angle $\sigma+\theta(\sigma)$ between the $x$-axis and the tangent to the free boundary.}
\label{Mestnikova-pic3}
\end{center}
\end{figure}

\medskip\noindent
\emph{Step~2.}
Now, we are able to determine the free boundary in the parametric form:
$$
\varGamma=\big\{(x,y)\in \mathbb{R}^2\;|\;
x=\widetilde{x}(\sigma),\; y=\widetilde{y}(\sigma),\; \sigma\in (0,\pi)\big\}
$$
with some functions $\widetilde{x}$ and $\widetilde{y}$ that will be
defined below. To this end, we recall that
$$
\varGamma=\big\{z\in \mathbb{C}\;|\; z=G\big(e^{i\sigma}\big),\; \sigma\in (0,\pi)\big\}
$$
and
\begin{equation*}
\frac{dG(e^{i\sigma})}{d\sigma}=
-  e^{-\tau(\sigma)} \, e^{i(\sigma+\theta(\sigma))}\,\cot\sigma.
\end{equation*}
Therefore, since the functions $\tau$ and $\theta$ are already known,
the functions $\widetilde{x}$ and $\widetilde{y}$ can be found from the
following relations:
\begin{subequations}\label{Mestnikova-2.28}
\begin{align}
\label{Mestnikova-2.28a}
& \frac{d\widetilde{x}(\sigma)}{d\sigma}=\MestnikovaRe\,\frac{dG\big(e^{i\sigma}\big)}{d\sigma}=
-e^{-\tau(\sigma)}\,\cos\big(\sigma+\theta(\sigma)\big)\,\cot\sigma,
\quad \widetilde{x}(\pi/2)=0,
\\
\label{Mestnikova-2.28b}
& \frac{d\widetilde{y}(\sigma)}{d\sigma}=\MestnikovaIm\,\frac{dG\big(e^{i\sigma}\big)}{d\sigma}=
-e^{-\tau(\sigma)}\,\sin\big(\sigma+\theta(\sigma)\big)\,\cot\sigma,
\quad \widetilde{y}(0)=1.
\end{align}
\end{subequations}
These formulas, in particular, imply that
\begin{equation*}
\frac{d\widetilde{y}(\sigma)}{d\sigma}
\Big(\frac{d\widetilde{x}(\sigma)}{d\sigma}\Big)^{-1}
=\tan\big(\sigma+\theta(\sigma)\big).
\end{equation*}
This means that $\sigma + \theta(\sigma)+\pi k$ is the
angle of inclination of $\varGamma$ at the point $\big(\widetilde{x}(\sigma),\widetilde{y}(\sigma)\big)$,
where $k$ is an integer number. The tangent vector $\boldsymbol\ell(\sigma)=\big(\widetilde{x}'(\sigma),\widetilde{y}'(\sigma)\big)$
is shown in Fig.~\ref{Mestnikova-pic3}. Since the inclination angle is equal to $\pi$ as $\sigma=0$,
we obtain that $k=1$ for $\sigma\in [0,\pi/2)$. The angle is equal to
$\pi$ as $\sigma=\pi$, therefore, $k=0$ for $\sigma\in (\pi/2,\pi]$.
The angle $\sigma+\theta(\sigma)$ does not have a discontinuity at $\sigma=\pi/2$ (see Fig.~\ref{Mestnikova-pic3}) and further,
in Section~\ref{Mestnikova-sec4.2}, it is called the inclination angle.

\medskip\noindent
\emph{Step~3.}
In order to determine the velocity field, we have to solve several standard problems. First of all, we find harmonic in the unit disk
$D_\circ^*=\{t\in \mathbb{C}\;|\; |t|< 1\}$ functions
$\hat{\tau}$ and $\hat{\theta}$ whose traces on the unit circle are
equal to $\tau$ and $\theta$, respectively.
Next, we define in $D_\circ^*$ the following functions:
\begin{gather*}
\varOmega(t)=\hat{\theta}(t)+i\hat{\tau}(t),\quad
 u(t)=-\frac{2}{\pi}\,\frac{1}{t}\, e^{-i\varOmega(t)},
\\
 w^*(t)=\frac{2}{\pi}\,\log\Big(\frac{t^2-1}{2t}\Big)-i,
\quad
 v^*(t)=\frac{dw^*(t)}{dt}=\frac{2}{\pi}\,\frac{1}{t}\,\frac{t^2+1}{t^2-1}
\end{gather*}
and find the mapping $z=G(t)$ as the solution of the following problem:
\begin{equation*}
\frac{dG(t)}{dt}=\frac{v^*(t)}{u(t)}=-\frac{t^2+1}{t^2-1}\,e^{i\varOmega(t)},
\quad G(0)=0.
\end{equation*}
The inverse mapping $F(z)$ can be determined as the solution of the problem:
$$
\frac{dF(z)}{dz}=\Big(G\,'(t)\Big)^{-1}\Big|_{t=F(z)}, \quad F(0)=0.
$$
Finally, the complex potential and the complex velocity in the plane of the variable $z$ can be found from the following relations:
$$
w(z)=w^*\big(F(z)\big),\quad v(z)=\frac{dw(z)}{dz}.
$$

Further in this paper, we restrict ourselves to solving the problem
of Step~1, i.e., to finding the functions $\tau$ and $\theta$.
The results of the numerical study of the problem are presented
in \cite{Mestnikova-MS}. Here, we investigate the existence and the uniqueness
of the solution (Section~\ref{Mestnikova-sec3}) as well as some of its properties (Section~\ref{Mestnikova-sec4}).

\section{Unique solvability of the problem}\label{Mestnikova-sec3}
This section is devoted to the investigation of the unique solvability of the problem \eqref{Mestnikova-2.25}, \eqref{Mestnikova-2.26} and \eqref{Mestnikova-2.27}.
Of course, one can study equations \eqref{Mestnikova-2.25} and \eqref{Mestnikova-2.26} directly, however, this system is rather difficult to deal with.
We have succeeded to find an equivalent formulation of the problem as a nonlinear operator equation in a Banach space and prove that it has a solution that is unique under a certain condition.

\subsection{Formulation of the problem in the form of an operator equation}\label{Mestnikova-sec3.1}
Equation \eqref{Mestnikova-2.26} includes a singular integral. It will be more convenient
to rewrite it in a form that contains an integral with an integrable kernel.
Notice that $\cot\sigma=\big(\log|\sin\sigma|\big)'$.
Therefore, using the integration by parts, the integral in equation \eqref{Mestnikova-2.26} can be represented as follows:
\begin{gather*}
\int_0^{\pi/2} \tau(s)\,\big(\cot(s-\sigma)-\cot(s+\sigma)\big)\,ds = \int_0^{\pi/2} \tau'(s)\, \big(\log|\sin(s+\sigma)| - \log|\sin(s-\sigma)| \big)\,ds
\\
=\int_0^{\pi/2} \tau'(s)\, \log \Big| \frac{\sin(s+\sigma)}{\sin(s-\sigma)} \Big|\,ds.
\end{gather*}
The integration by parts can be performed since we will find the solution with $\tau'\in L^2(0,\pi/2)$ which implies that $\tau$ will be H\"older continuous.
Thus, equation \eqref{Mestnikova-2.26} takes the form:
\begin{equation}\label{Mestnikova-3.1}
\theta(\sigma)=\frac{1}{\pi}\int_0^{\pi/2} K(s,\sigma)\,\tau'(s) \,ds,\quad \sigma\in [0,\pi/2],
\end{equation}
where
\begin{equation}\label{Mestnikova-3.2}
K(s,\sigma)=\log \Big| \frac{\sin(s+\sigma)}{\sin(s-\sigma)} \Big|.
\end{equation}
Notice that $K(s,0)=K(0,\sigma)=K(s,\pi/2)=K(\pi/2,\sigma)=0$ for $s,\sigma\in (0,\pi/2)$ and therefore
$\theta(0)=\theta(\pi/2)=0$. Thus, the boundary conditions
\eqref{Mestnikova-2.27} for $\theta$ follow from \eqref{Mestnikova-3.1}.

By making use of standard trigonometrical identities, the kernel $K$ can be
represented in another form:
\begin{equation} \label{Mestnikova-3.3}
K(s,\sigma)=\log\Big|\frac{\tan s+\tan\sigma}{\tan s-\tan\sigma}\Big|\, .
\end{equation}
In some situations, expression \eqref{Mestnikova-3.3} is more convenient than \eqref{Mestnikova-3.2}.
For instance, \eqref{Mestnikova-3.3} immediately implies that $K(s,\sigma)\geqslant 0$
for all values of $s$ and $\sigma$ in $(0,\pi/2)$ for which this kernel makes sense.

Let us introduce the following function
$$
\zeta(\sigma)=3\tau'(\sigma)
$$
and a nonlinear operator $\varPhi$ such that
\begin{equation*}
\varPhi(\zeta)(\sigma)=\frac{3}{\alpha\pi}\, \frac{\sin\big(\sigma+
\frac{1}{3}\,(H\zeta)(\sigma)\big)\,\cot\sigma}{\exp \int_0^\sigma \zeta\,ds},
\quad \sigma\in(0,\pi/2],
\end{equation*}
where
\begin{equation*}
(H\zeta)(\sigma)=\frac{1}{\pi}\int_{0}^{\pi/2} K(s,\sigma)\,\zeta(s) \, ds.
\end{equation*}
Equations \eqref{Mestnikova-2.25} and \eqref{Mestnikova-3.1} imply that the function $\zeta$
satisfies the following equation:
\begin{equation}\label{Mestnikova-3.4}
\zeta=\varPhi(\zeta).
\end{equation}
Really, as it follows from \eqref{Mestnikova-2.25},
$$
3\tau'(\sigma)=\frac{d}{d\sigma}\log \big(\alpha\pi + 3 \int_0^\sigma
\sin(s+\theta(s))\,\cot s\, ds\big).
$$
Integration of this equation from $0$ to an arbitrary $\sigma\in(0,\pi/2)$ and
the boundary condition \eqref{Mestnikova-2.27} imply that
$$
\alpha\pi + 3 \int_0^\sigma
\sin(s+\theta(s))\,\cot s\, ds=\alpha\pi\,\exp\big(\int_0^\sigma 3\tau'(s)\,ds\big).
$$
We have used here the H\"older continuity of $\theta$ which will be
justified in Section~\ref{Mestnikova-sec4}. This equality and equation \eqref{Mestnikova-2.25} immediately imply
\eqref{Mestnikova-3.4}.

Thus, the original problem is equivalent to the operator equation \eqref{Mestnikova-3.4}.
We will look for solutions of equation \eqref{Mestnikova-3.4} in the Hilbert space $L^2(0,\pi/2)$ equipped with the standard norm:
$$
\|\zeta\|=\Big(\int_0^{\pi/2}|\zeta(\sigma)|^2\,d\sigma\Big)^{1/2}.
$$
Further, for brevity, we will denote $L^2(0,\pi/2)$ by $L^2$.

\subsection{Auxiliary results}\label{Mestnikova-sec3.2}
At first, we recall some facts from the theory of Fourier series.
It is well known that the trigonometric systems
$\{\cos 2k\sigma\}_{k=0}^\infty$ and $\{\sin 2k\sigma\}_{k=1}^\infty$
are complete and orthogonal in $L^2$.
Therefore, for any function $f\in L^2$ and for almost all $\sigma\in [0,\pi/2]$,
the following representations hold true:
$$
f(\sigma)=\frac{a_0}{2}+\sum_{k=1}^\infty a_k \cos 2k\sigma,\quad
f(\sigma)=\sum_{k=1}^\infty b_k \sin 2k\sigma,
$$
where
$$
a_k=\frac{4}{\pi}\int_0^{\pi/2} f(s)\,\cos 2ks\,ds,\quad k=0,1,2,\ldots,
$$
$$
b_k=\frac{4}{\pi}\int_0^{\pi/2} f(s)\,\sin 2ks\,ds,\quad k=1,2,\ldots.
$$
For each of these expansions, we have the Bessel inequality and the Parseval equality:
$$
\frac{a_0^2}{2}+\sum_{k=1}^n a_k^2\leqslant \frac{4}{\pi}\, \|f\|^2,\quad
\sum_{k=1}^n b_k^2\leqslant \frac{4}{\pi}\, \|f\|^2,\quad n\in \mathbb{N},
$$
$$
\frac{a_0^2}{2}+\sum_{k=1}^\infty a_k^2=\frac{4}{\pi}\, \|f\|^2,\quad
\sum_{k=1}^\infty b_k^2= \frac{4}{\pi}\, \|f\|^2.
$$

For the kernel $K$, the following representation is well known:
$$
K(\sigma,s) = 2 \sum_{k = 1}^{\infty} \frac{\sin 2k\sigma \, \sin 2ks}{k}.
$$
It is true for all points, where the kernel $K(\sigma,s)$ is well defined, that is, for $\sigma + s$ and $\sigma - s$ that are different from $\pi m $ with an integer $m$.
 This representation is used in many works (see, for example, \cite{Mestnikova-PT} or \cite{Mestnikova-F}), however, we have
found its proof only in the book \cite{Mestnikova-N1}. This is a fairly rare book, for this reason, we give the proof in Appendix.
Notice that our proof differs from that in \cite{Mestnikova-N1}.  This representation will be useful in proving some estimates.

First, we study some properties of the operator $H$.
\begin{lemma}\label{Mestnikova-t3.1}
The linear operator $H:L^2\to L^2$ is bounded and
$$
\|Hw\|\leqslant \frac{1}{2}\,\|w\|
$$
for all $w\in L^2$.
\end{lemma}
\emph{Proof.}\hspace{3mm}
The boundedness of the operator $H$ in $L^2$ follows from the fact that
$$
\int_0^{\pi/2}\int_0^{\pi/2} \big|K(s,\sigma)\big|^2\, d\sigma ds<\infty.
$$
Let us prove the estimate in the assertion of the lemma.
Since $w\in L^2$, we have the following relation
\begin{equation}\label{Mestnikova-3.5}
w(\sigma)=\sum_{k=1}^{\infty} c_k \sin 2k\sigma
\end{equation}
that holds true for almost all $\sigma\in[0,\pi/2]$.
Here, $c_k$ are the Fourier coefficients. If $w_n$ is the $n$-th partial sum of this series, then
$$
(Hw_n)(\sigma)=
\sum_{k=1}^n \frac{c_k}{\pi}\int_{0}^{\pi/2} K(s,\sigma)\,\sin 2ks \, ds
=\sum_{k=1}^n \frac{c_k}{2k}\sin 2k\sigma .
$$
The Parseval equality and the Bessel inequality  imply that
$$
\frac{4}{\pi}\,\|Hw_n\|^2=\sum_{k=1}^n\frac{c_k^2}{4k^2}\leqslant
\frac{1}{4}\sum_{k=1}^n c_k^2\leqslant \frac{1}{\pi}\,\|w\|^2.
$$
The required estimate follows from this inequality, since $Hw_n\to Hw$ in $L^2$ as $n\to\infty$.
\hspace*{\fill}$\square$

\begin{lemma}\label{Mestnikova-t3.2}
If $w\in L^2$ and $u=Hw$, then $u'\in L^2$ and $\|u'\|= \|w\|$.
\end{lemma}
\emph{Proof.}\hspace{3mm}
We use the notations from the proof of Lemma~\ref{Mestnikova-t3.1}.
For every $w\in L^2$, the Fourier expansion \eqref{Mestnikova-3.5} takes place.
Denote again by $w_n$ the partial sum of this series and by $u_n$ the
function $Hw_n$. Due to Lemma~\ref{Mestnikova-t3.1}, $u_n\to u$ in $L^2$ as $n\to\infty$.
It is not difficult to see that
$$
u'_n(\sigma)=\sum_{k=1}^{n}c_k\,\cos 2k\sigma\quad\text{and}\quad
u'_n(\sigma)-u'_m(\sigma)=\sum_{k=m+1}^{n}c_k\,\cos 2k\sigma
$$
for $m<n$. As a consequence of the second equality, we find that
$$
\|u'_n-u'_m\|^2=\frac{\pi}{4}\sum_{k=m+1}^{n}c_k^2.
$$
Since the Fourier series in \eqref{Mestnikova-3.5} converges in $L^2$, the Cauchy criterion implies that the sequence $\{u_n'\}$ converges in $L^2$. Obviously, $u'$ is the limit.
The passage to the limit as $n\to\infty$ in the following equality
$$
\|u_n'\|^2=\frac{\pi}{4}\sum_{k=1}^{n}c_k^2
$$
enables us to conclude that $\|u'\|^2= \|w\|^2$ due to the Parseval equality.
\hspace*{\fill}$\square$

As a consequence of this lemma, we prove the following assertion.
\begin{lemma}\label{Mestnikova-t3.3}
If $w\in L^2$ and $u=Hw$, then $u\in C^{1/2}[0,\pi/2]$,
$$
|u(\sigma)-u(s)|\leqslant \|w\|\,|\sigma-s|^{1/2},\quad
|u(\sigma)|\leqslant \|w\|\,|\sigma|^{1/2},\quad\text{and}\quad
|u(\sigma)|\leqslant \|w\|\,|\pi/2-\sigma|^{1/2}
$$
for all $\sigma, s\in [0,\pi/2]$.
\end{lemma}
\emph{Proof.}\hspace{3mm}
Lemma~\ref{Mestnikova-t3.2} and the H\"older inequality imply that
$$
|u(\sigma)-u(s)|=\Big|\int_s^\sigma u'(t)\,dt\Big|\leqslant
\|u'\|\,|\sigma-s|^{1/2}\leqslant \|w\|\,|\sigma-s|^{1/2}.
$$
The second and the third inequalities in the assertion of the lemma
immediately follow from the first one since $u(0)=u(\pi/2)=0$ due to the properties of the kernel $K$.
\hspace*{\fill}$\square$

In the proof of the following lemma, we use some ideas form \cite{Mestnikova-F}.
\begin{lemma}\label{Mestnikova-t3.4}
If $w\in L^2$, then the functions
$$
f(\sigma)=\frac{(Hw)(\sigma)}{\sin\sigma}, \quad
h(\sigma)=\frac{\int_0^\sigma w(s)\,d s}{\sin\sigma},\quad
\text{and}\quad
g(\sigma)=  \sin\big(\sigma+\frac{1}{3} (Hw)(\sigma)\big)\,\cot\sigma
$$
belong to the space $L^2$ and the following estimates hold true:
\begin{gather}
\label{Mestnikova-3.6}
\|f\| \leqslant 2\|w\|,
\\
\label{Mestnikova-3.7}
\|h\| \leqslant 2\|w\|,
\\
\label{Mestnikova-3.8}
\|g\| \leqslant  1 + \frac{2}{3}\,\|w\|.
\end{gather}
\end{lemma}
\emph{Proof.}\hspace{3mm}
Let
$$
u_n(\sigma)=(Hw_n)(\sigma),\quad f_n(\sigma)=\frac{u_n(\sigma)}{\sin \sigma},\quad
\phi_{nm}(\sigma)=u_n(\sigma)-u_m(\sigma) \;\text{with}\; m<n,
$$
where $w_n$ is the $n$-th partial sum of the Fourier series \eqref{Mestnikova-3.5}.

At first, we show that the sequence $\{f_n\}$ converges in $L^2$.
To this end, we employ the Cauchy criterion. Let $\varepsilon$ be a small positive number.
Notice that
$$
\frac{1}{\sin^2 (s+\varepsilon)}= - \frac{d}{ds}\cot(s+\varepsilon)
$$
and, as it follows from Lemma~\ref{Mestnikova-t3.3}, $\phi_{nm}(0)=\phi_{nm}(\pi/2)=0$.
Making use of the integration by parts and the H\"older inequality,
we obtain the following estimate:
\begin{multline*}
\int_{0}^{\pi/2} \Big(\frac{\phi_{nm}(s)}{\sin(s+\varepsilon)} \Big)^2 \, ds =
\int_{0}^{\pi/2} 2\, \phi_{nm}(s)\,\phi'_{nm}(s)\, \cot(s+\varepsilon) \, ds \leqslant
\\
\leqslant \Big(\int_{0}^{\pi/2}\Big(\frac{\phi_{nm}(s)}{\sin(s+\varepsilon)}\Big)^2 \, ds
\Big)^{1/2}
\, \Big(\int_{0}^{\pi/2}\big(2\,\phi'_{nm}(s)\,\cos(s+\varepsilon) \big)^2 \, ds\Big)^{1/2}.
\end{multline*}
Since $\varepsilon$ is arbitrary, we have
$$
\int_{0}^{\pi/2} \Big(\frac{\phi_{nm}(s)}{\sin s} \Big)^2 \, ds\leqslant
4\int_{0}^{\pi/2}\big(\phi'_{nm}(s)\big)^2 \, ds=4\|\phi'_{nm}\|^2.
$$
Therefore, $\|f_n-f_m\|\leqslant 2\|\phi'_{nm}\|$.

Let us recall that $u'_n(\sigma)=\sum_{k=1}^{n}c_k\,\cos 2k\sigma$ and
$\phi'_{nm}=u'_n-u'_m=\sum_{k=m+1}^{n}c_k\,\cos 2k\sigma$.
Therefore,
\begin{equation*}
\|f_n-f_m\|^2\leqslant 4\|\phi'_{nm}\|^2= \pi\sum_{k=m+1}^{n} c_{k}^{2}.
\end{equation*}
Since $w\in L^2$, the series $\sum_{k=1}^{\infty} c_{k}^{2}$ converges.
The Cauchy criterion and the last estimate imply that the sequence  $\{f_n\}$ converges in $L^2$.
Moreover,
\begin{equation}\label{Mestnikova-3.9}
\|f_n\|^2\leqslant \pi\sum_{k=1}^{n} c_{k}^{2}\leqslant 4\|w\|^2
\end{equation}
for every $n\in \mathbb{N}$. Since $w_n\to w$ in $L^2$, Lemma~\ref{Mestnikova-t3.3} implies that $u_n=Hw_n\to Hw$ in $C^{1/2}[0,\pi/2]$. This means that $f_n(s)\to f(s)$ as $n\to\infty$ for every $s\in (0,\pi/2]$.
Due to \eqref{Mestnikova-3.9} and Fatou's lemma,  we conclude that $f\in L^2$ and $\|f\|\leqslant 2\|w\|$.
The first part of the lemma is proven.

The assertion concerning the function $h$ and inequality \eqref{Mestnikova-3.7} can be proven in exactly the same way.
In order to prove \eqref{Mestnikova-3.8},  notice that $|\sin x|\leqslant |x|$ for all $x\in \mathbb{R}$. Therefore,
$$
|g(\sigma)|\leqslant \big|\sigma+\frac{1}{3}\, (Hw)(\sigma)\big|\,\cot\sigma
\leqslant
\sigma\,\cot\sigma+\frac{1}{3}\,|f(\sigma)|\,\cos\sigma
$$
for all $\sigma\in (0,\pi/2]$. Since
$$
\int_0^{\pi/2} (\sigma\cot\sigma)^2\,d\sigma <1
$$
and $|\cos\sigma|\leqslant 1$, we have:
$$
\|g\| \leqslant 1 + \frac{1}{3} \,\|f(\sigma)\| \leqslant 1 + \frac{2}{3}\,\|w\|.
$$
The lemma is proven.
\hspace*{\fill}$\square$

\begin{lemma}\label{Mestnikova-t3.5}
For every $w\in L^1(0,\pi/2)$, define the function $W$ as follows:
$$
W(\sigma)=\int_{0}^{\sigma}w(s)\, ds.
$$
If $w\in L^2$, then $W \in C^{0, 1/2}[0,\pi/2]$ and
$|W(\sigma)|\leqslant \|w \| \sigma^{1/2}$ for all $\sigma\in [0,\pi/2]$.
\end{lemma}
\emph{Proof.}\hspace{3mm}
Since $W'(\sigma)=w(\sigma)$ for almost all $\sigma\in [0,\pi/2]$,
$$
|W(\sigma)-W(s)|=\Big|\int_s^\sigma w(t)\,dt\Big|\leqslant \|w\|\,|\sigma-s|^{1/2}.
$$
The inequality in the assertion of the lemma follows from the fact that $W(0)=0$.
\hspace*{\fill}$\square$

\begin{lemma}\label{Mestnikova-t3.6}
If $w_1$ and $w_2$ are non-negative functions from $L^2$, then
$$
\|\varPhi(w_2)-\varPhi(w_1)\|\leqslant \frac{8}{\alpha\pi}\,\|w_2-w_1\|.
$$
\end{lemma}
\emph{Proof.}\hspace{3mm}
Let $v_\lambda=w_1 + \lambda(w_2-w_1)$ for every $\lambda\in[0,1]$. Then
$$
\frac{dv_\lambda}{d\lambda}=w_2-w_1
$$
and
$$
\varPhi(w_2)-\varPhi(w_1)=\int_{0}^{1} \frac{d}{d\lambda} \varPhi(v_\lambda)d\lambda.
$$
It is not difficult to calculate that
\begin{multline*}
\frac{d}{d\lambda} \varPhi(v_\lambda)=\frac{3\cot\sigma}{\alpha\pi}\, \frac{d}{d\lambda} \frac{\sin(\sigma+\frac{1}{3}Hv_\lambda)}{\exp \int_{0}^{\sigma}v_\lambda \, ds} =
\\
=\frac{3\cot\sigma}{\alpha\pi}\,\Big(\frac{\cos(\sigma+\frac{1}{3}Hv_\lambda)\frac{1}{3}H(w_2-w_1)}{\exp \int_{0}^{\sigma}v_\lambda \, ds} -
\frac{\sin(\sigma+\frac{1}{3}Hv_\lambda)\int_{0}^{\sigma}(w_2-w_1)\, ds}{\exp \int_{0}^{\sigma}v_\lambda \, ds}\Big).
\end{multline*}
Taking into account the non-negativity of the function $v_\lambda$, we find that
$$
\Big\|\frac{d}{d\lambda} \varPhi(v_\lambda)\Big\|\leqslant  \frac{1}{\alpha\pi} \Big \|\frac{H(w_2-w_1)}{\sin\sigma} \Big \| + \frac{3}{\alpha\pi} \Big \|\frac{\int_{0}^{\sigma}(w_2-w_1)\, ds }{\sin\sigma}\Big \|.
$$
This inequality and Lemma~\ref{Mestnikova-t3.4} imply the following estimate:
$$
\Big\|\frac{d}{d\lambda} \varPhi(v_\lambda)\Big\|\leqslant
 \frac{2}{\alpha\pi}\, \|w_2-w_1\| + \frac{6}{\alpha\pi}\,\|w_2-w_1\| = \frac{8}{\alpha\pi}\, \|w_2-w_1\|,
$$
Now, the assertion of the lemma is a direct consequence of the inequality:
$$
\|\varPhi(w_2)-\varPhi(w_1)\|\leqslant \int_{0}^{1}
\Big \|\frac{d}{d\lambda}\varPhi(\upsilon_\lambda) \Big \| \, d\lambda.
$$
\hspace*{\fill}$\square$

\subsection{Unique solvability of the operator equation}\label{Mestnikova-sec3.3}
Now, we are able to prove  the existence and the uniqueness of the solution
of the operator equation  \eqref{Mestnikova-3.4}. At first, we investigate its solvability.
\begin{theorem}\label{Mestnikova-t3.7}
For every $\alpha\geqslant \frac{2}{\pi}\big(1+\frac{1}{3\sqrt{\pi}}\big)\approx 0.76$, there exists a non-negative function
$\zeta\in L^2$ that satisfies the operator equation \eqref{Mestnikova-3.4}. This function satisfies the estimates:
\begin{equation}\label{Mestnikova-3.10}
0\leqslant\sigma + \frac{1}{3}\, (H\zeta)(\sigma)\leqslant \pi \quad\text{for almost all}\quad\sigma\in (0,\pi/2),\quad
\|\zeta\|\leqslant \frac{9\sqrt{\pi}}{2}.
\end{equation}
\end{theorem}
\emph{Proof.}\hspace{3mm}
At first, we consider a sequence of the approximate problems
\begin{equation}\label{Mestnikova-3.11}
\zeta=\varPhi_k(\zeta),\quad k\in\mathbb{N},
\end{equation}
where
$$
\varPhi_k(\zeta)(\sigma)=\frac{3}{\alpha\pi}\, \frac{\sin\big(\sigma+ \frac{1}{3}\,
(H\zeta)(\sigma)\big)\,q_k(\sigma)}{\exp \int_0^\sigma \zeta\,ds},
\quad \sigma\in[0,\pi/2],
$$
and
$$
q_k(\sigma)=\begin{cases}
\displaystyle\cot \frac{1}{k}, & \sigma\in[0,1/k],
\\
\cot\sigma, & \sigma\in [1/k,\pi/2].
\end{cases}
$$
In order to prove the solvability of \eqref{Mestnikova-3.11}, we employ the Schauder fixed point theorem.
Let us define the set
$$
B_R^+=\{\zeta\in L^2\,| \,\|\zeta\|\leqslant R,\; \zeta\geqslant 0\}
$$
that is convex and closed in $L^2$. Our first goal is to establish the existence of  $R>0$
such that $\varPhi_k(B_R^+)\subset B_R^+$ for all $k\in \mathbb{N}$.

Let $R$ be a positive number and $\zeta\in B_R^+$. Then
$\exp\int_0^\sigma\zeta(s)\,ds \geqslant 1$ for all $\sigma\in [0,\pi/2]$.
Since $q_k(\sigma)\leqslant\cot\sigma$ for $k\in \mathbb{N}$ and $\sigma\in (0,\pi/2]$, Lemma~\ref{Mestnikova-t3.4} implies that
$$
\|\varPhi_k(\zeta)\|\leqslant \frac{3}{\alpha\pi}\,\big(1+\frac{2}{3}\,\|\zeta\|\big)
\leqslant \frac{3}{\alpha\pi}\,\Big(1+\frac{2}{3}\,R\Big).
$$
The inequality $\|\varPhi_k(\zeta)\|\leqslant R$ holds true, if $R$ is such that
$$
\frac{3}{\alpha\pi}\,\Big(1+\frac{2}{3}\,R\Big)\leqslant R.
$$
Thus, if
\begin{equation}\label{Mestnikova-3.12}
R\geqslant \frac{3}{\alpha\pi-2},
\end{equation}
then $\zeta\in B_R^+$ implies that $\|\varPhi_k(\zeta)\|\leqslant R$ for all $k\in \mathbb{N}$.

Next, we need to establish the non-negativity of the function $\varPhi_k(\zeta)$.
Since $q_k(\sigma)\geqslant 0$ for $\sigma\in[0,\pi/2]$, the non-negativity of the function $\varPhi_k(\zeta)$
is a consequence of the following inequality:
\begin{equation}\label{Mestnikova-3.13}
0\leqslant \sigma + \frac{1}{3}\, (H\zeta)(\sigma)\leqslant \pi.
\end{equation}
The left inequality holds true for a non-negative function $\zeta$, since
$K\geqslant 0$ and $H$ is a positive operator. As it follows from Lemma~\ref{Mestnikova-t3.3}, the right inequality
is valid whenever
$$
\|\zeta\|\,\gamma(\sigma)\leqslant 3(\pi-\sigma),
$$
where $\gamma(\sigma)=\min\{\sigma^{1/2},(\pi/2-\sigma)^{1/2}\}$. Therefore,
the function $\varPhi_k(\zeta)$ is non-negative, if $\zeta\in  B_R^+$ with
\begin{equation}\label{Mestnikova-3.14}
R\leqslant \min_{\sigma\in [0,\pi/2]}\frac{3(\pi-\sigma)}{\gamma(\sigma)}=
\frac{9\sqrt{\pi}}{2}.
\end{equation}

If $R_*$ is a number that satisfies inequalities \eqref{Mestnikova-3.12} and \eqref{Mestnikova-3.14},
then $\varPhi_k(B_{R^+_*})\subset B_{R^+_*}$ for all $k\in \mathbb{N}$.
Such a number $R_*$ exists, if
\begin{equation*}
\alpha\geqslant \frac{2}{\pi}\,\Big(1+\frac{1}{3\sqrt{\pi}}\Big)\approx 0.76.
\end{equation*}

Since $q_k(\sigma)\leqslant\cot\sigma$ for $k\in \mathbb{N}$ and $\sigma\in (0,\pi/2]$,
the estimate obtained in Lemma~\ref{Mestnikova-t3.6} holds true for the operators $\varPhi_k$, $k\in \mathbb{N}$.
This means that the operator $\varPhi_k$ is continuous in $L^2$ on non-negative functions for every
$k\in \mathbb{N}$. Let us prove its compactness. If $M$ is a bounded set in $B_{R^+_*}$, then
Lemmas~\ref{Mestnikova-t3.3} and \ref{Mestnikova-t3.5} imply that the function $\varPhi_k(\zeta)$ is H\"older continuous
on $[0,\pi/2]$ with the exponent $1/2$ whenever $\zeta\in M$. Therefore,
$\varPhi_k(M)$ is a bounded set in $C^{1/2}[0,\pi/2]$  for every $k\in \mathbb{N}$. Notice, however, that $\mathrm{diam}\,\varPhi_k(M)$
in $C^{1/2}[0,\pi/2]$, generally speaking, depends on $k$.
Since $C^{1/2}[0,\pi/2]$ is compactly embedded in $L^2$, the set $\varPhi_k(M)$ is compact in $L^2$ and,
as a consequence, the operator $\varPhi_k$ is compact on $B_{R^+_*}$ for every $k\in \mathbb{N}$ .

Thus,  $\varPhi_k$ satisfies all the conditions of the Schauder fixed point theorem for every  $k\in \mathbb{N}$ and
there exists a function $\zeta_k\in B_{R^+_*}\subset L^2$ that is a solution of equation \eqref{Mestnikova-3.11}.

The sequence $\{\zeta_k\}$ belongs to $B_{R^+_*}$ which is a bounded set in $L^2$. This means that
there exists a subsequence denoted again by $\{\zeta_k\}$ that weakly converges in $L^2$ to a function
 $\zeta\in B_{R^+_*}$. As follows from Lemmas~\ref{Mestnikova-t3.3} and \ref{Mestnikova-t3.5}, the sequences $\{H\zeta_k\}$ and
$\{\int_0^\sigma\zeta_k(s)\,ds\}$ are bounded in $C^{1/2}[0,\pi/2]$. Therefore, up to a subsequence,
$H\zeta_k\to H\zeta$ and $\int_0^\sigma\zeta_k(s)\,ds\to \int_0^\sigma\zeta(s)\,ds$ in $C[0,\pi/2]$
due to the compact embedding of $C^{1/2}[0,\pi/2]$ into $C[0,\pi/2]$. However, the sequence $\{\zeta_k\}$
does not converge in this space since $\zeta_k$ satisfies equation \eqref{Mestnikova-3.11} which includes the function $q_k$ and
the sequence $\{q_k\}$ does not converge in $C[0,\pi/2]$. Nevertheless,  $\{q_k\}$ converges to $\cot\sigma$ in
$C[\delta,\pi/2]$ for every $\delta>0$. Therefore,  $\{\zeta_k\}$ converges to $\zeta$ in $C[\delta,\pi/2]$ for every
$\delta>0$ and $\zeta(\sigma)= \varPhi(\zeta)(\sigma)$ for all $\sigma>\delta$.
Due to the arbitrariness of $\delta$, the function $\zeta$ is a solution of equation \eqref{Mestnikova-3.4}.

The estimates in \eqref{Mestnikova-3.10} follow from \eqref{Mestnikova-3.13} and \eqref{Mestnikova-3.14} respectively.
The theorem is proven.
\hspace*{\fill}$\square$

Our next goal is to prove the following uniqueness theorem.
\begin{theorem}\label{Mestnikova-t3.8}
For every $\alpha> \frac{8}{\pi}\approx 2.55$, there exists a unique non-negative function
$\zeta\in L^2$ that satisfies the operator equation \eqref{Mestnikova-3.4}.
\end{theorem}
\emph{Proof.}\hspace{3mm}
Let us denote by $\zeta_*$ the solution of \eqref{Mestnikova-3.4} whose existence is guaranteed by the previous theorem.
Since this solution is a non-negative function in $L^2$ and $\zeta_*=\varPhi(\zeta_*)$, Lemma~\ref{Mestnikova-t3.6}
implies that
$$
\|\varPhi(\zeta)-\zeta_*\|=\|\varPhi(\zeta)-\varPhi(\zeta_*)\|\leqslant \frac{8}{\alpha\pi}\,\|\zeta-\zeta_*\|
$$
for an arbitrary non-negative function $\zeta\in L^2$. If, in addition, $\zeta=\varPhi(\zeta)$, then
$$
\|\zeta-\zeta_*\|=\|\varPhi(\zeta)-\zeta_*\|\leqslant \frac{8}{\alpha\pi}\,\|\zeta-\zeta_*\|.
$$
For every $\alpha>\frac{8}{\pi}$, this inequality means that $\zeta=\zeta_*$ and, as a consequence, that
equation \eqref{Mestnikova-3.4} has only one non-negative solution in $L^2$.
\hspace*{\fill}$\square$

\section{Solution of the original problem}\label{Mestnikova-sec4}

As noted in Section~\ref{Mestnikova-sec2.3}, the solution of the original problem can be constructed,
if we know the functions $\tau$ and $\theta$. Suppose that $\zeta_*\in L^2$ is the non-negative
solution of the operator equation \eqref{Mestnikova-3.4} whose existence was proven in the previous section.
Let us define the functions $\tau$ and $\theta$ as follows:
\begin{equation}\label{Mestnikova-4.1}
\theta(\sigma)=\frac{1}{3}\,(H\zeta_*)(\sigma),\quad \tau'(\sigma)=\frac{1}{3}\,\zeta_*(\sigma),
\quad \sigma\in (0,\pi/2).
\end{equation}
These formulas entirely correspond to the notations introduced in Section~\ref{Mestnikova-sec3.1}.
Notice that we have a differential equation for the function $\tau$, for this reason, we need
the boundary condition $\tau(0)=\log\pi/2$ that is taken from \eqref{Mestnikova-2.27}.
The functions $\tau$ and $\theta$ defined above satisfy
equations \eqref{Mestnikova-2.25} and \eqref{Mestnikova-2.26}. According to Lemmas~\ref{Mestnikova-t3.3} and \ref{Mestnikova-t3.5}, these
functions are in $C^{0,1/2}[0,\pi/2]$.

\subsection{Smoothness of the solution}\label{Mestnikova-sec4.1}
The functions $\tau$ and $\theta$ have, in fact, a higher smoothness than is stated above.
At first, we prove two auxiliary lemmas.

\begin{lemma}\label{Mestnikova-t4.1}
If $w$ is a non-negative function from $L^2$ and $v(\sigma)=\sin\sigma\,\varPhi(w)(\sigma)$ for $\sigma\in(0,\pi/2)$, then
$v'\in L^2$ and $\|v'\|\leqslant \frac{3}{\alpha\pi}\,\big(\sqrt{\pi} + \frac{4}{3}\,\|w\|\big)$.
\end{lemma}
\emph{Proof.}\hspace{3mm}
It is not difficult to see that
\begin{multline*}
v'(\sigma)=\frac{3}{\alpha\pi}\,\frac{1}{\exp\int_0^\sigma w(s)\,ds}\,\Big(
\cos(\sigma + \frac{1}{3}Hw)\,\big(1+\frac{1}{3}(Hw)'\big)\,\cos\sigma
\\
-\sin(\sigma + \frac{1}{3}Hw)\,\big(\sin\sigma +w(\sigma)\,\cos\sigma\big)\Big).
\end{multline*}
Since $w$ is assumed to be non-negative, Lemma~\ref{Mestnikova-t3.2} implies that $v'\in L^2$.
The estimate follows from the fact that $\|(Hw)'\|=\|w\|$.
\hspace*{\fill}$\square$

\begin{lemma}\label{Mestnikova-t4.2}
If $w\in L^2$, $(w\,\sin \sigma)'\in L^2$, and $w(\pi/2)=0$,  then $\big((Hw)'\sin\sigma\big)'\in L^2$
and $\|\big((Hw)'\sin\sigma\big)'\|\leqslant\|(w\,\sin\sigma)'\|$.
\end{lemma}
\emph{Proof.}\hspace{3mm}
Let $w(\sigma)=\sum_{k=1}^{\infty} c_k \sin 2k\sigma$
and $w_n$ be the $n$-th partial sum of this series.
Using the well-known trigonometric formulas, it is not difficult to find that
$$
w_n(\sigma) \sin\sigma =\frac{1}{2} \,\sum_{k=1}^{n} c_k \big(\cos (2k-1)\sigma - \cos (2k+1)\sigma\big)
=\sum_{k=1}^{n+1}\gamma_k \cos(2k-1)\sigma,
$$
where
$$
\gamma_1=\frac{c_1}{2}, \quad \gamma_k=\frac{1}{2}\,(c_k-c_{k-1}), \quad
k=2,\ldots,n,\quad \gamma_{n+1}=-\frac{c_n}{2}.
$$
The function systems $\{\cos(2k-1)\sigma\}_{k\in \mathbb{N}}$ and $\{\sin(2k-1)\sigma\}_{k\in \mathbb{N}}$
are, generally speaking, not complete in $L^2$. At least, we cannot provide the prove of this fact. If it were so, the inequality
in the assertion of the lemma might be replaced by the equality. However, these systems are orthogonal in $L^2$ and
$$
\int_0^{\pi/2} \cos^2(2k-1)\sigma\,d\sigma=\int_0^{\pi/2} \sin^2(2k-1)\sigma\,d\sigma=\frac{\pi}{4},
\quad k\in \mathbb{N}.
$$
Since $w_n\to w$ in $L^2$ as $n\to\infty$, we have the following equalities:
$$
w(\sigma) \sin\sigma =\sum_{k=1}^{\infty}\gamma_k \cos(2k-1)\sigma,
\quad
\sum_{k=1}^\infty\gamma_k^2=\frac{4}{\pi}\,\|w\sin\sigma\|^2.
$$

Let us calculate the Fourier coefficients of the function $(w\sin\sigma)'\in L^2$
with respect to the system $\{\sin(2k-1)\sigma\}_{k\in \mathbb{N}}$:
$$
\int_0^{\pi/2} \big(w(\sigma) \sin\sigma\big)'\, \sin(2k-1)\sigma\,d\sigma=
-(2k-1) \int_0^{\pi/2} w(\sigma) \sin\sigma\, \cos(2k-1)\sigma\,d\sigma=
-\frac{\pi}{4}\,\gamma_k\,(2k-1).
$$
Due to the Bessel inequality
\begin{equation}\label{Mestnikova-4.2}
\sum_{k=1}^{\infty}\gamma_k^2 (2k-1)^2\leqslant \frac{4}{\pi}\|(w\sin\sigma)'\|^2.
\end{equation}

Exactly as in the proof of Lemma~\ref{Mestnikova-t3.1}, we find that
$$
(Hw_n)' (\sigma)\sin\sigma = \frac{1}{2}\,\sum_{k=1}^{n} c_k \big(\sin(2k+1)\sigma - \sin(2k-1)\sigma\big)
=-\sum_{k=1}^{n+1}\gamma_k \sin(2k-1)\sigma
$$
and, as a consequence, that
$$
\big((Hw_n)'(\sigma)\sin \sigma\big)'=-\sum_{k=1}^{n+1}\gamma_k(2k-1) \cos(2k-1)\sigma.
$$
This  equality and \eqref{Mestnikova-4.2} give the following estimate
$$
\|\big((Hw_n)'\sin \sigma\big)'\|\leqslant \|(w\sin\sigma)'\|, \quad n\in \mathbb{N}.
$$
Since $\big((Hw_n)'\sin \sigma\big)'$ converges weakly
to $\big((Hw)'\sin \sigma\big)'$ in $L^2$ as $n\to\infty$, this estimate implies
the assertion of the lemma.
\hspace*{\fill}$\square$

Now, we apply these lemmas to proving the differentiability of the functions $\tau$ and $\theta$.
\begin{theorem}\label{Mestnikova-t4.3}
Let $\tau$ and $\theta$ be the functions defined in \eqref{Mestnikova-4.1}. Then\\
1. $(\tau'\sin\sigma)'\in L^2$ and $\|(\tau'\sin\sigma)'\|\leqslant 7/(\alpha\sqrt{\pi})$, $\tau'(\pi/2)=0$, $\tau\in C^{1,1/2}[\delta,\pi/2]$ for every $\delta\in(0,\pi/2)$.\\
2. $(\theta'\sin\sigma)'\in L^2$ and $\|(\theta'\sin\sigma)'\|\leqslant 7/(\alpha\sqrt{\pi})$,  $\theta\in C^{1,1/2}[\delta,\pi/2]$ for every $\delta\in(0,\pi/2)$.
\end{theorem}
\emph{Proof.}\hspace{3mm}
At first, we note that  $\zeta_*(\pi/2)=0$, since $\zeta_*$ satisfies \eqref{Mestnikova-3.4} and $\cot\pi/2=0$.
 Lemma~\ref{Mestnikova-t4.1} implies that the function
$\big(\zeta_*(\sigma)\sin\sigma\big)'$ is in $L^2$ and $\|(\zeta_*\sin\sigma)'\|\leqslant \frac{3}{\alpha\pi}\,\big(\sqrt{\pi} + \frac{4}{3}\,\|\zeta_*\|\big)$. Due to the second estimate in \eqref{Mestnikova-3.10},
$\|(\zeta_*\sin\sigma)'\|\leqslant \frac{21}{\alpha\sqrt{\pi}}$. Further, as a consequence of
Lemma~\ref{Mestnikova-t4.2}, we conclude that $\big((H\zeta_*)'\sin\sigma\big)'\in L^2$ and
$\|((H\zeta_*)'\sin\sigma)'\|\leqslant \|(\zeta_*\sin\sigma)'\|\leqslant \frac{21}{\alpha\sqrt{\pi}}$.
The assertion of the theorem follows from \eqref{Mestnikova-4.1} and the embedding theorems.
\hspace*{\fill}$\square$

\medskip\noindent
\textbf{Remark.}
Theorem~\ref{Mestnikova-t4.3} states the smoothness of the functions $\tau$ and $\theta$ on $[\delta,\pi/2]$ for all
$\delta\in (0,\pi/2)$. However, the symmetry conditions \eqref{Mestnikova-2.24a} enable us to define
these functions also on $[0,\pi]$. Thus, $\tau,\theta\in C^{1,1/2}[\delta,\pi-\delta]$ for every $\delta\in(0,\pi/2)$.
\hfill\textbullet

\medskip
The smoothness of the functions $\tau$ and $\theta$ can be further improved, namely, they are analytic.
In order to prove this fact, we employ the classical result by H.~Lewy \cite{Mestnikova-L}
which is formulated below as a lemma. Notice that, in this lemma and in the proof of the subsequent theorem, the cartesian variables
$x$ and $y$ are different from the physical variables used in the beginning of our paper. These symbols were used in \cite{Mestnikova-L}.

\begin{lemma}\label{Mestnikova-t4.4}
Let $U(x,y)$ be harmonic near the origin in $y<0$ and $V(x,y)$ be the conjugate harmonic of $U$.
Assume that $U$, $V$, $\partial U/\partial x$ exist and are continuous in the semi-neighborhood
$y\leqslant 0$ of the origin. If the boundary values on $y=0$ satisfy a relation
$$
\frac{\partial U}{\partial y}= A\big(x, U, V, \partial U/\partial x\big)
$$
in which $A$ is an analytic function of all four arguments for all values occurring, then
$U(x,y)$ and $V(x,y)$ are analytically extensible across $y=0$.
\end{lemma}

\begin{theorem}\label{Mestnikova-t4.5}
The functions $\tau$ and $\theta$ are analytic on $(0,\pi)$.
\end{theorem}
\emph{Proof.}\hspace{3mm}
The functions $\tau$ and $\theta$ are the traces of the conjugate harmonic in the unit disk $D_\circ^*$
functions $\hat{\tau}$ and $\hat{\theta}$. Besides that, these functions
are continuously differentiable on $(0,\pi)$ and satisfy equation \eqref{Mestnikova-2.15}. The Cauchy --- Riemann
conditions written in the polar coordinates $(\varrho,\sigma)$ together with \eqref{Mestnikova-2.15} imply that
\begin{equation}\label{Mestnikova-4.3}
\frac{\partial\hat{\theta}(t)}{\partial\varrho}=A_0(\sigma,\hat{\theta},\hat{\tau})
\end{equation}
for $\varrho=|t|=1$ and $\sigma\in (0,\pi)$, where $A_0(\sigma,\hat{\theta},\hat{\tau})=
\frac{\pi^2}{8\alpha}\,e^{-3\hat{\tau}(t)} \sin(\sigma+\hat{\theta}(t))\,\cot\sigma$ is an analytic
function of its arguments.

Let us take an arbitrary $\sigma_0\in (0,\pi)$. There exists a conformal linear fractional transformation
$(x,y) \mapsto t=f(x,y)$ that maps the lower half-plane $y<0$ to the unit disk $D_\circ^*$ and such that
$f(0)=e^{i\sigma_0}$. If we introduce $U(x,y)=\hat{\theta}(f(x,y))$ and $V(x,y)=\hat{\tau}(f(x,y))$,
then, in the new variables, equation \eqref{Mestnikova-4.3} takes the following form:
$$
\frac{\partial U}{\partial y}= A\big(x, U, V, \partial U/\partial x\big),
$$
with an analytic function $A$. Lemma~\ref{Mestnikova-t4.4} implies that the functions $U(x,0)$ and $V(x,0)$ are analytic.
Since the inverse mapping $f^{-1}$ is conformal, the functions $\theta$ and $\tau$ are analytic in a
neighborhood of the point $\sigma_0$. The assertion of the theorem follows from the arbitrariness
of $\sigma_0$.
\hspace*{\fill}$\square$

\subsection{The form of the free boundary}\label{Mestnikova-sec4.2}
\begin{figure}[h]
\begin{center}
\includegraphics[width=0.25\textwidth]{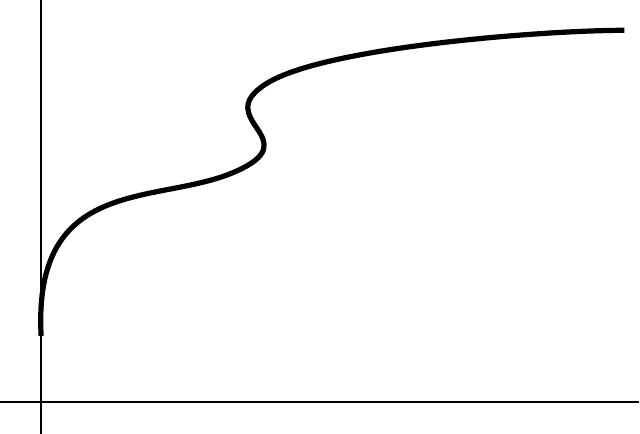}
\hspace{2cm}
\includegraphics[width=0.46\textwidth]{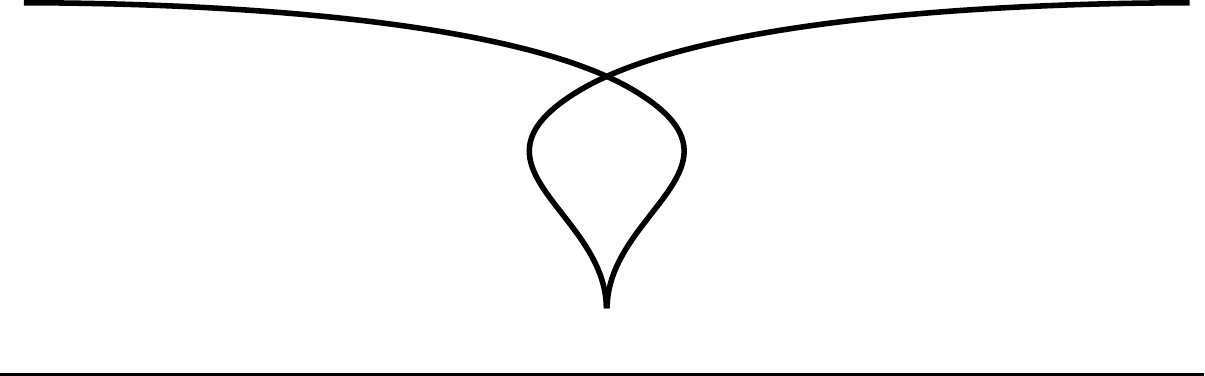}
\\
\parbox[t]{0.25\textwidth}{\caption{\small The possible form of the free boundary with  $\sigma_0+\theta(\sigma_0)\geqslant \pi/2$
for some $\sigma_0\in (0,\pi/2)$.}
\label{Mestnikova-pic4}}
\hspace{2cm}
\parbox[t]{0.46\textwidth}{\caption{\small The non-injectivity of the mappings $G$ and $F$ in the case of overturning at $\sigma$ close to $\pi/2$.}
\label{Mestnikova-pic5}}
\end{center}
\end{figure}

We intend to investigate the inclination angle $\sigma+\theta(\sigma)$ of the free boundary.
We denote this angle by $\eta$:
$$
\eta(\sigma)=\sigma+\theta(\sigma)=\sigma +\frac{1}{3}\,(H\zeta_*)(\sigma), \quad \sigma\in[0,\pi/2].
$$
Due to the properties of the operator $H$, $\theta(0)=\theta(\pi/2)=0$, which implies that
$\eta(0)=0$ and $\eta(\pi/2)=\pi/2$.
The solution $\zeta_*$ of equation \eqref{Mestnikova-3.4} whose existence is stated in
Theorem~\ref{Mestnikova-t3.7} satisfies \eqref{Mestnikova-3.10}. This means that $0\leqslant \eta(\sigma)\leqslant \pi$
for all $\sigma\in [0,\pi/2]$.
As it follows from \eqref{Mestnikova-2.28b}, the function $\widetilde{y}(\sigma)$
decreases monotonically on the interval $(0,\pi/2)$, however, the function $\widetilde{x}(\sigma)$
can behave in such a way that overturning occurs. This situation is shown in Fig.~\ref{Mestnikova-pic4}.
Moreover, the overturning at $\sigma$ close to $\pi/2$ implies that the mappings $G$ and $F$ are
not one-to-one (see Fig.~\ref{Mestnikova-pic5}). This means that such a solution of the operator equation \eqref{Mestnikova-3.4},
if it exists, does not give the solution of the original problem.

Here, we show that such a situation is impossible,
namely, the inclination angle $\eta(\sigma)$ of the free boundary
is less than $\pi/2$, if the reduced Froude number $\alpha$ is not too small.
This question is often encountered in problems with a free boundary (see, e.g., \cite{Mestnikova-KN,Mestnikova-PT}).
In order to obtain the estimate of such a kind, the maximum principle for a subharmonic function is usually employed.
We failed to implement this approach due to the presence of the singular sink that causes the unboundedness of
corresponding functions. For this reason, we suggest another proof of this fact based on estimates of the solution.
Of course, these estimates are not optimal and can be improved.

For brevity, we introduce the following notation:
$$
w(\sigma) = \sin\eta(\sigma)\,\cot\sigma,\quad\sigma\in(0,\pi/2].
$$
It is not difficult to see that $w(\sigma)>0$ for $\sigma\in (0,\pi/2)$, $w(\pi/2)= 0$, and, as it follows from Lemma~\ref{Mestnikova-t3.4}, $w\in L^2$.

\begin{lemma}\label{Mestnikova-t4.6}
The function $(w(\sigma)\sin\sigma)'$ is in $L^2$ and
$$
\|(w\sin\sigma)'\|\leqslant  \sqrt{\pi} +\frac{1}{3}\,\|\zeta_*\|.
$$
\end{lemma}
\emph{Proof.}\hspace{3mm}
Since
$$
(w(\sigma)\sin\sigma)'=(\sin\eta(\sigma)\,\cos\sigma)'=\eta'(\sigma)\,\cos\eta(\sigma)\,\cos\sigma
-\sin\eta(\sigma)\,\sin\sigma
$$
and $\eta'=1+(H\zeta_*)'/3$, we have the following estimate:
$$
|(w(\sigma)\sin\sigma)'|\leqslant |\cos\sigma|+\frac{1}{3}\,|(H\zeta_*)'(\sigma)|+|\sin\sigma|, \quad \sigma\in[0,\pi/2].
$$
Therefore, due to Lemma~\ref{Mestnikova-t3.2}, we have:
$$
\|(w\sin\sigma)'\|\leqslant \sqrt{\pi} +\frac{1}{3}\,\|(H\zeta_*)'\|\leqslant \sqrt{\pi} +\frac{1}{3}\,\|\zeta_*\|.
$$
\hspace*{\fill}$\square$

\begin{theorem}\label{Mestnikova-t4.7}
If $\alpha>5/\sqrt{8}\approx 1.768$, then $\eta(\sigma)<\pi/2$ for all $\sigma\in (0,\pi/2)$.
\end{theorem}
\emph{Proof.}\hspace{3mm}
Since $\zeta_*$ is a positive solution of \eqref{Mestnikova-3.4} and
$H$ is a positive operator,
\begin{equation}\label{Mestnikova-4.4}
\eta(\sigma)=\sigma +\frac{1}{\alpha\pi}\,H\Big(\frac{\sin\eta\,\cot\sigma}{\exp\int_0^\sigma\zeta_*\,ds}\Big)
\leqslant \sigma +\frac{1}{\alpha\pi}\,H(\sin\eta\,\cot\sigma)= \sigma +\frac{1}{\alpha\pi}\,Hw.
\end{equation}
Thus, $\eta(\sigma)<\pi/2$ if
\begin{equation}\label{Mestnikova-4.5}
(Hw)(\sigma)<\alpha\pi\,(\pi/2-\sigma),\quad \sigma\in(0,\pi/2).
\end{equation}
Our goal is to prove this estimate.

Due to Lemma~\ref{Mestnikova-t4.2},
$$
\|((Hw)'\sin\sigma)'\|\leqslant \|(w\sin\sigma)'\|.
$$
Since $(Hw)'(\sigma)\sin\sigma=0$ for $\sigma=0$, we have the following estimate:
$$
|(Hw)'(\sigma)\sin\sigma|=\Big|\int_0^\sigma ((Hw)'(s)\sin s)'\,ds\Big|\leqslant
\sigma^{1/2}\,\|((Hw)'\sin\sigma)'\|\leqslant \sigma^{1/2}\,\|(w\sin\sigma)'\|
$$
for all $\sigma\in [0,\pi/2]$. Together with Lemma~\ref{Mestnikova-t4.6} and the second estimate in \eqref{Mestnikova-3.10}, this implies that
$$
|(Hw)'(\sigma)|\leqslant \|(w\sin\sigma)'\|\, \frac{\sigma^{1/2}}{\sin\sigma}\leqslant A\, \phi(\sigma),\quad \sigma\in (0,\pi/2],
$$
where
$$
A=\frac{5}{2}\,\sqrt{\pi}\quad\text{and}\quad \phi(\sigma)= \frac{\sigma^{1/2}}{\sin\sigma}.
$$
Notice that $\phi(\pi/2)=\phi(\pi/4)=\sqrt{\pi/2}$ and $\phi(\sigma)< \sqrt{\pi/2}$ for $\sigma\in (\pi/4,\pi/2)$.
Therefore, since $(Hw)(\pi/2)=0$,
\begin{equation}\label{Mestnikova-4.6}
|(Hw)(\sigma)|=\Big|\int_{\sigma}^{\pi/2} (Hw)'(s)\,ds\Big|
\leqslant A\int_{\sigma}^{\pi/2}\phi(s)\,ds
\leqslant A\,\sqrt{\frac{\pi}{2}}\,\Big(\frac{\pi}{2}-\sigma\Big), \quad \sigma\in [\pi/4,\pi/2].
\end{equation}
Thus, estimate \eqref{Mestnikova-4.5} holds true for $\sigma\in [\pi/4,\pi/2]$, if
\begin{equation}\label{Mestnikova-4.7}
\alpha>A/\sqrt{2\pi}=5/\sqrt{8}\approx 1.768.
\end{equation}

In order to obtain a similar estimate for $\sigma\in (0,\pi/4)$, we apply other arguments.
Due to the definition of the function $w$ and \eqref{Mestnikova-4.4}, we have:
$$
w(\sigma)\leqslant \eta(\sigma)\,\cot\sigma\leqslant \sigma\cot\sigma + \frac{1}{\alpha\pi}\,\cos\sigma\,\frac{(Hw)(\sigma)}{\sin\sigma}.
$$
As it follows from Lemma~\ref{Mestnikova-t3.4},
$$
\|w\|\leqslant B +\frac{2}{\alpha\pi}\,\|w\|,
$$
where $B$ is the $L^2$-norm of the function $\sigma\cot\sigma$. Notice that
$B^2=\pi \log 2 -\pi^3/24\approx 0.89<1$.
Therefore,
$$
\|w\|\leqslant \frac{\alpha\pi}{\alpha\pi - 2}.
$$
Lemma~\ref{Mestnikova-t3.3} implies that
$$
|(Hw)(\sigma)|\leqslant \|w\|\,\sigma^{1/2}\leqslant \frac{\alpha\pi}{\alpha\pi - 2}\,\sigma^{1/2}.
$$
Thus, estimate \eqref{Mestnikova-4.5} is true for $\sigma\in(0,\pi/4)$, if
$$
\frac{\alpha\pi}{\alpha\pi - 2}\,\sqrt{\frac{\pi}{4}}\leqslant \alpha\pi\,\Big(\frac{\pi}{2}-\frac{\pi}{4}\Big)
$$
which is equivalent to the following inequality:
\begin{equation}\label{Mestnikova-4.8}
\alpha\geqslant \frac{2}{\pi^{3/2}}+\frac{2}{\pi}\approx 0.996.
\end{equation}

Estimates \eqref{Mestnikova-4.7} and \eqref{Mestnikova-4.8} imply that \eqref{Mestnikova-4.5} holds true for $\alpha$ that satisfies
the inequality \eqref{Mestnikova-4.7}. The theorem is proven.
\hspace*{\fill}$\square$

Let us make a remark to the result obtained. One can see a significant difference between inequalities
\eqref{Mestnikova-4.7} and \eqref{Mestnikova-4.8} deduced for $\sigma\in [\pi/4,\pi/2]$ and $\sigma\in (0,\pi/4)$,
respectively. It may seem that we can improve the assertion of the theorem, if we split the interval $(0,\pi/2)$ by a number
that differs from $\pi/4$. However, our method of proof does not enable us to considerably
improve estimate \eqref{Mestnikova-4.6} and, as a consequence, \eqref{Mestnikova-4.7}. The point is that the function $\phi$ in \eqref{Mestnikova-4.6} is almost constant on $[\pi/4,\pi/2]$ and does not differ much from $\sqrt{\pi/2}$.
Moreover, $\phi(\pi/2)=\sqrt{\pi/2}$.

The estimates obtained in the proof of the previous theorem make it possible to determine the asymptotics
of the inclination angle as $\sigma\to\pi/2$, i.e., in the neighborhood of the cusp point.
\begin{theorem}\label{Mestnikova-t4.8}
If $\alpha>\alpha_0=5/\sqrt{8}\approx 1.768$, i.e., $\alpha$ satisfies the same inequality as in Theorem~\ref{Mestnikova-t4.7}, then
$$
\eta(\sigma)=\pi/2 - \beta\,(\pi/2-\sigma) +o(\pi/2-\sigma)\quad
\text{as}\quad \sigma\to\pi/2^{-},
$$
where $\beta\in [1-\alpha_0/\alpha, 1]$ is a constant.
\end{theorem}
\emph{Proof.}\hspace{3mm}
As it follows from Theorem~\ref{Mestnikova-t4.3} (and from Theorem~\ref{Mestnikova-t4.4}), the function $\eta$ is
differentiable and $\beta$ is its derivative at the point $\sigma=\pi/2$.
The main assertion of the theorem is that $\beta\ne 0$. Thus, we have to estimate
$\beta=\eta'(\pi/2)$. At first, we estimate $\beta$ from below. Due to \eqref{Mestnikova-4.4} and \eqref{Mestnikova-4.6},
we have:
\begin{multline*}
\beta=\lim_{\sigma\nearrow\pi/2}\frac{\eta(\sigma) -\eta(\pi/2)}{\sigma -\pi/2}\geqslant
\lim_{\sigma\nearrow\pi/2}\frac{\sigma +(\alpha\pi)^{-1}\,(Hw)(\sigma)-\pi/2}{\sigma -\pi/2}
\\
=1+\frac{1}{\alpha\pi}\lim_{\sigma\nearrow\pi/2}\frac{(Hw)(\sigma)}{\sigma -\pi/2}
\geqslant 1-\frac{1}{\alpha\pi}\,A\sqrt{\frac{\pi}{2}}=1-\frac{\alpha_0}{\alpha}.
\end{multline*}
In order to prove the upper bound for $\beta$, we notice that $\eta(\sigma)=\sigma+\theta(\sigma)$
and $\theta(\sigma)\geqslant 0$ for $\sigma\in [0,\pi/2]$. Therefore,
$$
\beta=\lim_{\sigma\nearrow\pi/2}\frac{\sigma+\theta(\sigma) -\pi/2}{\sigma -\pi/2}
=1 + \lim_{\sigma\nearrow\pi/2}\frac{\theta(\sigma)}{\sigma -\pi/2}\leqslant 1.
$$
\hspace*{\fill}$\square$

Notice that the estimates in the previous proof imply that $-\alpha_0/\alpha\leqslant \theta'(\pi/2)\leqslant 0$.
This means that $\theta'(\pi/2)$ tends to $0$ as $\alpha\to\infty$.

It would be also interesting to study the asymptotics of the free boundary in the neighborhood of the cusp point
in the original variables $x$ and $y$.  We prove the following assertion.
\begin{theorem}\label{Mestnikova-t4.9}
Let the condition of Theorem~\ref{Mestnikova-t4.8} be satisfied, i.e., $\alpha>\alpha_0=5/\sqrt{8}\approx 1.768$.
If the free boundary $\varGamma$ is described
by the equation $y=q(x)$ for $x\geqslant 0$ with some function $q$, then
\begin{equation}\label{Mestnikova-4.9}
q(x)= y_0+ a\, x^{2/3}+o(x^{2/3})\quad\text{as}\quad x\to 0,
\end{equation}
where $y_0=\widetilde{y}(\pi/2)$ is the $y$-coordinate of the cusp point,
$a=3^{2/3}\sqrt{c_0}/(2\beta^{2/3})$, $c_0=e^{-\tau(\pi/2)}$,  and $\beta$ is defined in Theorem~\ref{Mestnikova-t4.8}.
\end{theorem}
\emph{Proof.}\hspace{3mm}
The free boundary is described parametrically by the functions $\widetilde{x}(\sigma)$ and $\widetilde{y}(\sigma)$ which satisfy \eqref{Mestnikova-2.28}.
The right-hand sides in \eqref{Mestnikova-2.28} are analytic functions of $\sigma$. If we leave in their power expansions only the leading-order terms,
we obtain that
\begin{subequations}\label{Mestnikova-4.10}
\begin{align}
& \frac{d\widetilde{x}(\sigma)}{d\sigma}=-e^{-\tau(\sigma)}\,\cos\eta(\sigma)\,\cot\sigma
= -c_0\beta(\pi/2-\sigma)^2 +o\big((\pi/2-\sigma)^2\big),
\\
& \frac{d\widetilde{y}(\sigma)}{d\sigma}=-e^{-\tau(\sigma)}\,\sin\eta(\sigma)\,\cot\sigma
= -c_0(\pi/2-\sigma)+o\big((\pi/2-\sigma)^2\big)
\end{align}
\end{subequations}
as $\sigma\to \pi/2$, where $c_0=e^{-\tau(\pi/2)}$ and $\beta$ is defined in Theorem~\ref{Mestnikova-t4.8}.
Here, we have used Theorem~\ref{Mestnikova-t4.8} and the following relations:
\begin{align*}
& \cot\sigma=\pi/2-\sigma +O\big((\pi/2-\sigma)^3\big),
\\
& \cos\eta(\sigma)=\sin\big(\beta(\pi/2-\sigma)+o(\pi/2-\sigma)\big)=\beta(\pi/2-\sigma)+o(\pi/2-\sigma),
\\
& \sin\eta(\sigma)=\cos\big(\beta(\pi/2-\sigma)+o(\pi/2-\sigma)\big)=1+ O\big((\pi/2-\sigma)^2\big),
\\
& e^{-\tau(\sigma)}=c_0  + O\big((\pi/2-\sigma)^2\big)
\end{align*}
as $\sigma\to \pi/2$. Since $\widetilde{x}(\pi/2)=0$ and $\widetilde{y}(\pi/2)=y_0$, equations \eqref{Mestnikova-4.10}
imply that
\begin{subequations}\label{Mestnikova-4.11}
\begin{align}
\label{Mestnikova-4.11a}
& \widetilde{x}(\sigma)= \frac{c_0\beta}{3}\,(\pi/2-\sigma)^3 +\omega_x\big((\pi/2-\sigma)^3\big),
\\
\label{Mestnikova-4.11b}
& \widetilde{y}(\sigma)= y_0+\frac{c_0}{2}\,(\pi/2-\sigma)^2+\omega_y\big((\pi/2-\sigma)^3\big),
\end{align}
\end{subequations}
where $\omega_x(s)=o(s)$ and $\omega_y(s)=o(s)$ as $s\to 0$. Notice that the functions
$\omega_x$ and $\omega_y$ are differentiable and $\omega_x'(0)=\omega_y'(0)=0$.
Due to the inverse function theorem, there exists a differentiable function $s=s(x)$ such that
$$
\frac{c_0\beta}{3}\,s(x) +\omega_x\big(s(x)\big)=x
$$
for $x$ close to $0$. This means that equation \eqref{Mestnikova-4.11a} is satisfied with $(\pi/2-\sigma)^3=s(x)$.
Besides that, since $s(0)=0$ and $s'(0)=3/(c_0\beta)$,
$$
s(x)=\frac{3}{c_0\beta}\, x + o(x)\quad\text{as}\quad x\to 0.
$$
Substitution of these relations into \eqref{Mestnikova-4.11b} gives:
\begin{multline*}
q(x)=y_0+\frac{c_0}{2}\,\Big(\frac{3}{c_0\beta}\, x + o(x)\Big)^{2/3} +
\omega_y\big(\frac{3}{c_0\beta}\, x + o(x)\big)
\\
=y_0 + \frac{3^{2/3}}{2}\,\frac{\sqrt{c_0}}{\beta^{2/3}}\, x^{2/3} + o(x^{2/3})\quad\text{as}\quad x\to 0.
\end{multline*}
The theorem is proven.
\hspace*{\fill}$\square$

\appendix

\section*{Appendix}\label{Mestnikova-sec5}
\refstepcounter{section}
In Section~\ref{Mestnikova-sec3.2}, we used the well known representation of the kernel $K$ as the sum of a trigonometric series.
The proof of this representation is rather difficult to find and, for completeness, we give it below.

\begin{lemma}\label{Mestnikova-t5.1}
Assume that $\sigma, s\in \mathbb{R}$ and the numbers $\sigma+s$ and $\sigma-s$ are
not equal to $\pi m$, $m\in \mathbb{Z}$. If
$$
K(\sigma,s)=\ln \Big| \frac{\sin(\sigma+s)}{\sin(\sigma-s)} \Big|,
$$
then
$$
K(\sigma,s)=2 \sum_{k=1}^{\infty} \frac{\sin 2k\sigma \, \sin 2ks}{k}.
$$
\end{lemma}
\emph{Proof.}\hspace{3mm}
By making use of trigonometric identities, we find that
$$
2\sum_{k=1}^\infty \frac{\sin 2k\sigma \, \sin 2k s}{k} = \sum_{k=1}^\infty \frac{\cos 2k(\sigma -s)}{k} - \sum_{k=1}^\infty \frac{\cos 2k(\sigma + s)}{k}=A(\sigma -s)-A(\sigma + s),
$$
where
$$
A(\varphi)=\sum_{k=1}^\infty \frac{\cos 2k\varphi}{k}.
$$
According to the Dirichlet test, this series converges for all $\varphi\ne \pi m$, $m\in \mathbb{Z}$,
since
$$
\sum_{k=1}^n \cos 2k\varphi =\frac{\cos (n+1)\varphi \, \sin n\varphi}{\sin\varphi},
\quad n\in \mathbb{N}.
$$
Therefore, $A(\varphi)$ is well defined for $\varphi\ne \pi m$, $m\in \mathbb{Z}$.

Let us introduce the following function  $S$ of the complex variable $z$:
$$
S(z) = \sum_{k=1}^\infty \frac{z^k}{k}.
$$
It is not difficult to see that $A(\varphi)= \MestnikovaRe S(e^{2 i \varphi})$ for all $\varphi\ne \pi m$, $m\in \mathbb{Z}$.
At the same time,
$$
S'(z) = \sum_{k=1}^\infty z^{k-1} = \frac{1}{1-z}\quad\text{as}\quad |z|<1.
$$
This means that $S(z)=-\log (1-z)$ as $|z|<1$. The function $S$ is also well defined
everywhere on the circle $|z|=1$ except at the point $z=1$. Thus,
extending $S(z)$ by continuity, we obtain that
$$
A(\varphi)=-\ln 2|\sin \varphi|\quad\text{for all}\quad\varphi\ne \pi m,\; m\in \mathbb{Z},
$$
which immediately implies the assertion of the lemma.
\hspace*{\fill}$\square$

\end{document}